\documentclass [11pt]{article}

\newtheorem{theorem}{Theorem}
\newtheorem{definition}[theorem]{Definition}
\newtheorem{lemma}[theorem]{Lemma}
\newtheorem{proposition}[theorem]{Proposition}
\newtheorem{corollary}[theorem]{Corollary}
\newtheorem{question}[theorem]{Question}
\usepackage{amsmath}
\usepackage{amsfonts}
\usepackage{mathrsfs}
\usepackage{amssymb}
\usepackage{amscd}
\usepackage{scalerel}
\usepackage{graphicx}
\usepackage{url}
\usepackage[lined,algonl,boxed,norelsize]{algorithm2e}

\newcommand{\Ccal}{{\mathcal C}}

\newcommand{\Fcal}{{\mathcal F}}

\newcommand{\Ocal}{{\mathcal O}}

\newcommand{\Tcal}{{\mathcal T}}

\newcommand{\Div}{\operatorname{Div}}
\newcommand{\Prin}{\operatorname{Prin}}

\newcommand{\indeg}{\operatorname{indeg}}
\newcommand{\Jac}{\operatorname{Jac}}

\newcommand{\Pic}{\operatorname{Pic}}

\newcommand{\Lspan}{\operatorname{span}}

\newcommand{\sign}{\operatorname{sign}}

\newcommand{\la}{{\langle}}
\newcommand{\ra}{{\rangle}}

\DeclareRobustCommand{\rchi}{{\mathpalette\irchi\relax}}
\newcommand{\irchi}[2]{\raisebox{\depth}{$#1\chi$}}

\title{Geometric Bijections Between Spanning Trees and Break Divisors}

\author{Chi Ho Yuen}

\begin{document}
\maketitle

\begin{abstract}
The Jacobian group $\Jac(G)$ of a finite graph $G$ is a group whose cardinality is the number of spanning trees of $G$. $G$ also has a tropical Jacobian which has the structure of a real torus; using the notion of break divisors, An et al. obtained a polyhedral decomposition of the tropical Jacobian where vertices and cells correspond to elements of $\Jac(G)$ and spanning trees of $G$, respectively. We give a combinatorial description of bijections coming from this geometric setting. This provides a new geometric method for constructing bijections in combinatorics. We introduce a special class of geometric bijections that we call edge ordering maps, which have good algorithmic properties. Finally, we study the connection between our geometric bijections and the class of bijections introduced by Bernardi; in particular we prove a conjecture of Baker that planar Bernardi bijections are ``geometric''. We also give sharpened versions of results by Baker and Wang on Bernardi torsors.
\end{abstract}

\section{Introduction} \label{sec:Intro}

Tropical geometry is the study of algebraic geometry in terms of a ``combinatorial shadow''; for instance, the study of algebraic curves becomes a study of graphs and metric graphs. Many notions in classical algebraic geometry now have combinatorial counterparts, such as {\em divisors}, {\em Jacobians}, and the {\em Riemann-Roch theorem} for graphs and metric graphs \cite{BF_Metric}, \cite{BN_RR}, \cite{GK_RR}. Furthermore, tropical geometry also possesses some combinatorial tools that do not exist in the classical world. One example is the notion of {\em break divisors} introduced by Mikhalkin and Zharkov \cite{MZ_TropJac}. They proved that there is a continuous section to the map $\Div_+^g(\Gamma)\rightarrow\Pic^g(\Gamma)$ from degree $g$ effective divisors to the degree $g$ divisor classes on a metric graph $\Gamma$ of genus $g$; the set of break divisors is the image of this section. Such a section does not exist in the world of algebraic curves. A combinatorial proof of the Mikhalkin-Zharkov result was given by An, Baker, Kuperberg and Shokrieh \cite{ABKS}, who also proved a discrete version of the theorem for finite graphs, namely that {\em integral break divisors} form a set of representatives for $\Pic^g(G)$.\\

Break divisors are closely related to {\em spanning trees} of a graph and other combinatorial concepts. On finite graphs, the set of (integral) break divisors has the same cardinality as the set of spanning trees, and break divisors are equivalent to {\em indegree sequences of root-connected orientations} and Gioan's {\em cycle-cocycle reversal classes} \cite{Gioan_CCsys}. On metric graphs, the authors of \cite{ABKS} constructed a canonical {\em polyhedral decomposition} of $\Pic^g(\Gamma)$ where each full-dimensional cell canonically represents a spanning tree. Using this polyhedral decomposition, they gave a ``volume proof''\footnote{The concept of a ``volume proof'' of matrix--tree theorem was known in previous literature (e.g. \cite{Stanley_zono}), but their constructions are non-canonical as they need to fix extra data other than the graph itself.} of the classical Kirchhoff's matrix--tree theorem. Such a polyhedral decomposition is also related to tropical geometry topics including Brill-Noether theory and tropical theta divisors; we refer the reader to \cite[Section 5.4]{BJ_Survey} for details.\\

The authors of \cite{ABKS} observed that generic ``shiftings'' of their polyhedral decomposition induce bijections between integral break divisors and spanning trees. In this paper we give a combinatorial description of such bijections.  These ``geometric bijections'' are special cases of what we call {\em cycle orientation maps}.  We give a necessary and sufficient condition for a cycle orientation map to be geometric, and hence a sufficient condition for a cycle orientation map to be bijective. Our proof is geometric in nature, and gives a new ``non-combinatorial'' approach to proving bijectivity; indeed, we do not know a combinatorial proof of bijectivity for general geometric bijections. We also study algorithmic aspects of a particular class of cycle orientation maps that we call {\em edge ordering maps}.\\

Bernardi introduced a process on graphs using {\em ribbon structures} \cite{Bernardi_Process}, that is, combinatorial data to specify embeddings of graphs on orientable surfaces. Using his process, Bernardi obtained various bijections between certain classes of subgraphs and certain classes of (indegree sequences of) orientations, as well as a new characterization of the Tutte polynomial. In particular, Bernardi gave a bijection between spanning trees and indegree sequences of root-connected orientations. Baker observed that Bernardi's bijections map a cell of the aformentioned polyhedral decomposition to one of its vertices, which is also a feature of the bijections induced from geometric shiftings. On the other hand, treating those indegree sequences as elements of $\Pic^g(G)$, Bernardi's bijections induce simply transitive group actions (or {\em torsor structures}) of the Jacobian on the set of spanning trees. Motivated by a question of Ellenberg \cite{MathOverFlow}, \cite{AIM}, Baker and Wang proved that if the ribbon structure is planar, then all Bernardi torsors are isomorphic and they respect plane duality \cite{BW_Torsor}. Similar results were proven for {\em rotor routing} by Chan et al. \cite{CCG}, \cite{CGM} (see \cite{HLMPPW} for background on rotor routing). In fact, Baker and Wang showed that the torsors induced by Bernardi bijections are isomorphic to the torsors induced by rotor routing if the ribbon structure is planar. Based on these facts as well as some computational evidence, Baker asked whether planar Bernardi bijections come from the above geometric picture \cite[Remark 5.2]{BW_Torsor}. In this paper, we answer Baker's question in the affirmative, and we give alternative proofs and sharpenings of Baker and Wang's results.\\ 

The paper is structured as follows. In Section \ref{sec:Background}, we give some background for various topics we consider in this paper; \ref{sec:Intro_Graph} is fundamental for all of the paper, \ref{sec:Intro_Metric} and \ref{sec:Intro_Decomp} are used in Section \ref{sec:GB}, \ref{sec:Intro_GO} is the combinatorial background for Section \ref{sec:Bijection} and Section \ref{sec:Bernardi}, \ref{sec:Intro_Bernardi} and \ref{sec:Intro_Duality} contain material needed for Section \ref{sec:Bernardi}. In Section \ref{sec:GB}, we give our combinatorial description of geometric bijections. In Section \ref{sec:Bijection}, we focus on edge ordering maps, which have good combinatorial and algorithmic properties. In Section \ref{sec:Bernardi}, we study when a Bernardi bijection is geometric and give new proofs of some fundamental facts about Bernardi bijections. In Section \ref{sec:Problems}, we mention several open problems for future research. Section \ref{sec:Bijection} is purely combinatorial, and can be read independently of the geometric discussion in Section \ref{sec:GB}; Section \ref{sec:Bernardi} is almost independent of the previous two sections except definitions from Definition \ref{Def_COC} and Theorem \ref{Def_EOM}, and the statement of Theorem \ref{main_theorem}. Therefore readers interested in particular parts can read the corresponding background and jump to the individual sections directly.\\

{\bf Acknowledgement}: Many thanks to Matt Baker for proposing the conjecture relating geometric bijections and Bernardi bijections to the author, as well as for the helpful comments and suggestions throughout the author's research and writing process. The author also wants to thank Spencer Backman, Farbod Shokrieh, Olivier Bernardi, Greg Kuperberg and Lionel Levine for engaging conversations, Robin Thomas and Josephine Yu for checking some of the technical graph theory and polyhedral geometry lemmas in this paper, and Emma Cohen for drawing Figure \ref{Bernardi_OM_diag}. Finally the author thanks the anonymous referees for their detailed comments and suggestions.

\section{Background} \label{sec:Background}

\subsection{Graphs and Divisors} \label{sec:Intro_Graph}

Unless otherwise specified, all graphs in this paper are assumed to be finite and connected, possibly with parallel edges but without any loops; the term cycle will refer to a simple cycle. We use $n$ and $m$ to denote the number of vertices and edges of a graph, respectively. The {\em Laplacian matrix} $\Delta$ of a graph is the $n\times n$ matrix defined as $D-A$, where $D$ is a diagonal matrix with the $(i,i)$-entry being the degree of the vertex $v_i$, and $A$ is the adjacency matrix of $G$ with the $(i,j)$-entry being the number of edges between $v_i$ and $v_j$ for $i\neq j$ and 0 for $i=j$.\\

A {\em divisor} on a graph $G=(V,E)$ is a function $D:V\rightarrow\mathbb{Z}$. The set of all divisors on $G$ is denoted by $\Div(G)$; it has a natural group structure as the free abelian group generated by $V$. The {\em degree} $\deg(D)$ of a divisor $D$ is $\sum_{v\in V}D(v)$, and the set of all divisors of degree $k$ is denoted by $\Div^k(G)$. Identifying divisors with vectors in $\mathbb{Z}^n$, we say a divisor is {\em principal} if it is equal to $\Delta {\bf u}$ for some integral vector ${\bf u}$. The set of all principal divisors is denoted by $\Prin(G)$. $\Prin(G)$ is a subgroup of $\Div^0(G)$ and the quotient group $\Div^0(G)/\Prin(G)$ is the (degree 0) {\em Picard group} $\Pic^0(G)$ of $G$; the Picard group of a graph is also known as the Jacobian, the sandpile group or the critical group of the graph in other literature. More generally we say two divisors $D,D'$ are {\em linearly equivalent}, denoted $D\sim D'$, if $D-D'\in\Prin(G)$. We denote by $\Pic^k(G)$ the set $\Div^k(G)/\!\!\sim$. It is easy to see that $\Pic^k(G)$ and $\Pic^0(G)$ differ only by a degree $k$ translating element.\\

Given a divisor $D$, we often say there are $D(v)$ {\em chips} at vertex $v$, and we can construct a divisor on $G$ by adding or removing chips on vertices in various ways. The concept of linear equivalence can be captured by the {\em chip-firing game} on $G$, see \cite{BN_RR} for its usage in the Riemann-Roch theory on graphs and \cite{LP_AMS} for a brief overview of chip-firing in other parts of mathematics.\\

A {\em spanning tree} of $G$ is a connected spanning subgraph of $G$ with no cycles. Denote by $g:=m-n+1$ the {\em genus} (or cyclomatic number) of $G$, which is the number of edges outside any spanning tree of $G$\footnote{We clarify that $g$ is not the topological genus of $G$ as in the usual topological graph theory sense.}. It is known that the Picard group has the same cardinality as the set of spanning trees $S(G)$ of a graph \cite{BHN_Lattice, DM_ASM}. A degree $g$ divisor is called an {\em (integral) break divisor} if it can be obtained from the following procedure: choose a spanning tree $T$ of $G$, and for every edge $f\not\in T$, pick an orientation of $f$ and add a chip at the head of $f$. In general, different choices of spanning trees and orientations can produce the same break divisor. It is proven in \cite[Theorem 1.3]{ABKS} that break divisors form a set of representatives for the divisor classes in $\Pic^g(G)$. In particular, the number of break divisors is equal to the number of spanning trees of $G$.

\subsection{Metric Graphs} \label{sec:Intro_Metric}

A {\em metric graph} $\Gamma$ (or an abstract tropical curve) is a compact connected metric space such that every point has a neighborhood isometric to a star-shaped set. Every metric graph can be constructed in the following way: start with a weighted graph $(G=(V,E),w:E\rightarrow\mathbb{R}_{>0})$, associate with each edge $e$ a closed line segment $L_e$ of length $w(e)$, and identify the endpoints of the $L_e$'s according to the graph structure in the obvious way. A (weighted) graph $G$ that yields a metric graph $\Gamma$ is said to be a {\em model} of $\Gamma$. The {\em genus} of $\Gamma$ is the genus of any model of $\Gamma$. See \cite{BF_Metric} for details. For simplicity, we assume all graphs we consider have uniform edge weights 1, though most geometric propositions in this paper remain valid for the general case.\\

Metric graphs have an analogous theory of divisors as graphs; we refer the reader to \cite{ABKS} for details and only explain the most relevant notions here. A {\em divisor} on $\Gamma$ is an element of the free abelian group generated by $\Gamma$; equivalently, it is a function $D:\Gamma\rightarrow\mathbb{Z}$ of finite support. A divisor $D=(p_1)+(p_2)+\ldots+(p_g)$ is a {\em break divisor} if there is a model $G$ for $\Gamma$ and an edge $e_i$ of $G$ for each $p_i$ such that $p_i\in L_{e_i}\subset\Gamma$ (by the definition of $L_{e_i}$, $p_i$ can be on the endpoints of $e_i$), and $G-\{e_1,\ldots,e_g\}$ is a spanning tree of $G$. As in the case of graphs, break divisors on $\Gamma$ form a set of representatives for the equivalence classes in $\Pic^g(\Gamma)$, cf. \cite[Theorem 1.1]{ABKS}, \cite{MZ_TropJac}.\\

We say a divisor on a metric graph $\Gamma$ is {\em integral} with respect to a model $G$ if the support of $D$ consists of vertices of $G$. For any model $G$ of $\Gamma$, the break divisors of $G$ naturally correspond to the integral break divisors on $\Gamma$ with respect to $G$.

\subsection{Tropical Jacobians and Their Polyhedral Decompositions} \label{sec:Intro_Decomp}

The Picard group $\Pic^0(\Gamma)$, and hence any $\Pic^k(\Gamma)$, has a natural $g$-dimensional real torus structure. We give an informal description here; readers interested in a rigorous treatment can refer to \cite{BF_Metric}.\\

Fix a model $G$ for $\Gamma$, as well as an arbitrary orientation for the edges of $G$, which we will call the {\em reference orientation}. The {\em real edge space} $C_1(G,\mathbb{R})$ of $G$ is the $m$-dimensional vector space with $E$ as a basis over $\mathbb{R}$. This space is equipped with an {\em inner product} $\la \cdot,\cdot\ra$ that extends $\la e_i,e_j\ra=\delta_{i,j}w(e_i)$ bilinearly. Given a cycle $C$ (resp. a path $P$) of $G$, we can interpret $C$ as an element of $C_1(G,\mathbb{R})$: pick an orientation of $C$ and take the $\pm 1$-combination of edges in $C$, where the coefficient of $e\in C$ is 1 precisely when the reference orientation of $e$ agrees with its orientation induced from the chosen orientation of $C$. In general there are two possible choices of orientation for the cycle (resp. path), but for the rest of this paper either an arbitrary choice works, or we will explicitly describe which orientation to choose. The {\em real cycle space} $H_1(G,\mathbb{R})$ is the $g$-dimensional subspace of $C_1(G,\mathbb{R})$ spanned by all cycles of $G$; denote by $\pi:C_1(G,\mathbb{R})\rightarrow H_1(G,\mathbb{R})$ the orthogonal projection.\\

The projections of $k$ edges $f_1,\ldots,f_k\in G$ onto $H_1(G,\mathbb{R})$ are linearly independent if and only if $G-\{f_1,\ldots,f_k\}$ is connected, and in particular $\pi(e_1),\ldots,\pi(e_g)$ form a basis of $H_1(G,\mathbb{R})$ exactly when $e_1,\ldots,e_g$ are the edges outside some spanning tree $T$. Let $C_i^T$ denotes the unique cycle (known as the {\em fundamental cycle}) contained in $T+e_i$, then $C_1^T,\ldots,C_g^T$ is a basis of $H_1(G,\mathbb{R})$. The {\em integral cycle space} $H_1(G,\mathbb{Z})$ is the $g$-dimensional lattice in $H_1(G,\mathbb{R})$ consisting of integral combinations of cycles of $G$; it is generated by $C_1^T,\ldots,C_g^T$ for the set of fundamental cycles of any spanning tree \cite{Biggs_APG}. The {\em tropical Jacobian} $\Jac(\Gamma)$ is the $g$-dimensional real torus $H_1(G,\mathbb{R})/H_1(G,\mathbb{Z})$ with the induced inner product.\\

Now we define a bijection $\Phi$ between $\Pic^g(\Gamma)$ and $\Jac(\Gamma)$. Fix a vertex $q$ of $G$. For each divisor class in $\Pic^g(\Gamma)$, pick the unique break divisor $D=(p_1)+\ldots+(p_g)$ from it\footnote{In fact, any representative yields the same image.}. For each $p_i$, choose a path $\gamma_i$ from $q$ to $p_i$ and interpret $\gamma_i$ as an element of $C_1(G,\mathbb{R})$. Then the image of $[D]$ in $\Jac(\Gamma)$ is $[\pi(\gamma_1+\ldots+\gamma_g)]\in H_1(G,\mathbb{R})/H_1(G,\mathbb{Z})=\Jac(\Gamma)$. Different choices of $q$ produce the same $\Phi$ up to translations in the universal cover of $\Jac(\Gamma)$ only, so essentially $\Phi$ is independent of $q$, and from now on we will often abuse notations and identify $\Pic^g(\Gamma)$ and $\Jac(\Gamma)$ using $\Phi$.\\

\begin{figure}[h!]
\begin{center}
    \includegraphics[scale=0.4]{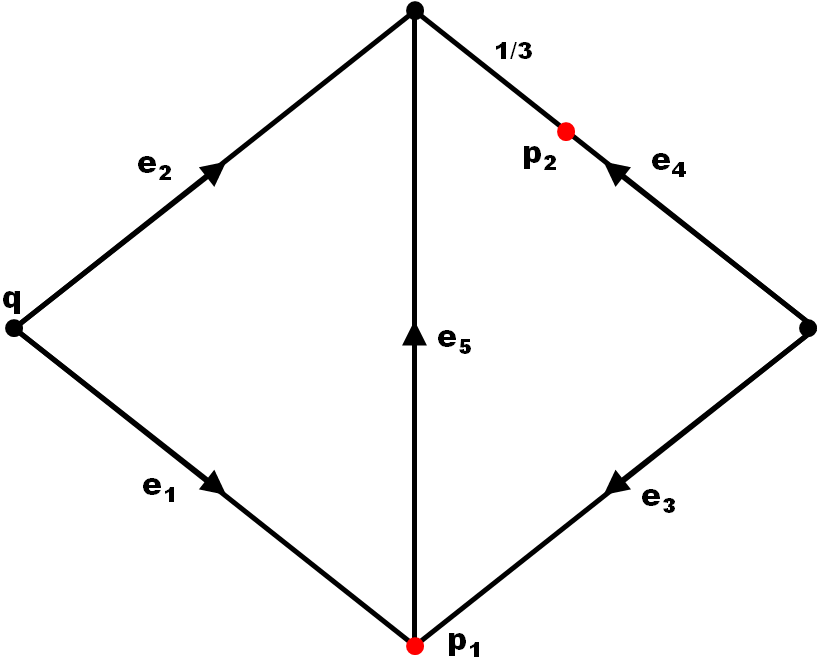}
\end{center}
  \caption{A fixed model for a metric graph $\Gamma$ and a break divisor $(p_1)+(p_2)$ on it.}
  \label{AJ_diag}
\end{figure}

\noindent{\bf Example.} In Figure \ref{AJ_diag}, take the spanning tree $\{e_1,e_3,e_5\}$. The cycle space is spanned by $C_1:=e_1-e_2+e_5,C_2:=e_3-e_4+e_5$. Take the paths $\gamma_1=e_1,\gamma_2=e_2-\frac{1}{3}e_4$ from $q$ to $p_1,p_2$. The image of $(p_1)+(p_2)$ in $\Jac(\Gamma)$ is $[\pi((e_1)+(e_2-\frac{1}{3}e_4))]=[\frac{-1}{24}C_1+\frac{1}{8}C_2]=[\frac{23}{24}C_1+\frac{1}{8}C_2]$.\\

It is interesting to ask how the image of $\Phi$ changes when a break divisor is perturbed by a small amount. Let $D=(p_1)+\ldots+(p_g)$ be a break divisor, choose a spanning tree $T$ with $e_1=\overrightarrow{u_1v_1},\ldots,e_g=\overrightarrow{u_gv_g}$ denoting the (unit length) edges not in $T$ in such a way that $p_i\in L_{e_i}$. Suppose that for some $j$ the segment $\overrightarrow{u_jp_j}$ is of length $\theta_j<1$, let $D'=(p_1)+\ldots+(p'_j)+\ldots+(p_g)$ where $p'_j$ is a point in $e_j$ such that $\overrightarrow{u_jp'_j}$ is of length $\theta_j<\theta'_j\leq 1$. Then $\Phi([D'])=\Phi([D])+[(\theta'_j-\theta_j)\pi(e_j)]$.\\

A generic break divisor with respect to a model $G$, that is a break divisor with each point $p_i$ of its support lying in the interior of some edge $e_i$, determines a unique spanning tree $T=G-\{e_1,\ldots,e_g\}$. The set of all generic break divisors that determine the same spanning tree $T$ forms a connected open set in $\Pic^g(\Gamma)$. More precisely, let $T$ be a spanning tree and let the edges not in $T$ be $e_1=\overrightarrow{u_1v_1},\ldots,e_g=\overrightarrow{u_gv_g}$. Then by definition $D=(u_1)+\ldots+(u_g)$ is a break divisor and the image of every generic break divisor compatible with $T$ is equal to $\Phi([D])+[\sum_{i=1}^g \theta_i\pi(e_i)]$ for some $0<\theta_i<1$. The closure $\Ccal_T$ of such subset is therefore the image of (a translation of) the parallelotope $\{\sum_{i=1}^g\theta_i\pi(e_i):0\leq\theta_i\leq 1\}\subset H_1(G,\mathbb{R})$ in $\Jac(\Gamma)$. We call $\Ccal_T$ the {\em cell corresponding to $T$}. The tropical Jacobian is the union of $\Ccal_T$ as $T$ runs through all spanning trees of $G$. As two distinct cells are disjoint except possibly at the boundary, they give a {\em polyhedral decomposition} of $\Pic^g(\Gamma)\cong\Jac(\Gamma)$. Also, by construction, the vertices of such decomposition are exactly the integral break divisors on $\Gamma$. For proof of these properties, we refer the reader to \cite[Section 3]{ABKS}.\\

\begin{figure}[h!]
\begin{center}
    \includegraphics[width=0.5\textwidth]{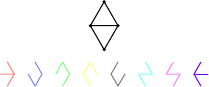}
\end{center}
  \vspace{+10pt}
\begin{center}
    \includegraphics[width=0.8\textwidth]{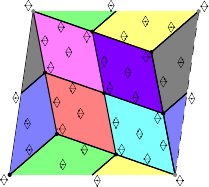}
\end{center}
  \caption{The polyhedral decomposition of $\Pic^2(\Gamma)$ (\cite[Figure 1]{ABKS}).}
  \label{Decomp_diag}
\end{figure}

{\bf Remark.} As an outcome of studying the combinatorial meaning of the polyhedral decomposition, we state a combinatorial interpretation of all faces of the polyhedral decomposition here, which might be useful for future work. We say a pair $(D',E')$, with $D'\in\Div^i(G),E'\subset E(G), |E'|=j$, is a {\em break $(i,j)$-configuration} if $G-E'$ is connected and $D'$ is a break divisor of $G-E'$. Thus, for examples, a $(g,0)$-configuration specifies a break divisor and a $(0,g)$-configuration specifies a spanning tree. For each $0\leq i\leq g$, there is a one-to-one correspondence between the $i$-dimensional faces and the break $(g-i,i)$-configurations, where a break configuration $(D',E')$ corresponds to the face consisting of break divisors of the form $D'+D''$ where each chip of $D''$ is in an edge of $E'$.

\subsection{Graph Orientations and Divisors} \label{sec:Intro_GO}

In the paper of An-Baker-Kuperberg-Shokrieh \cite{ABKS}, many key results related to break divisors were proven using the concept of {\em orientable divisors}, which was further developed into an algorithmic theory by Backman \cite{Bac_GO}. In this section, basic definitions and results are reviewed. We mostly work with full orientations on finite graphs, though the general theory includes partial orientations on metric graphs.

\begin{definition}
Let $\Ocal$ be an orientation of a graph $G$. The divisor $D_{\Ocal}$ is $\sum_{v\in V(G)}(\indeg(v)-1)(v)$, so $\deg(D_{\Ocal})=g-1$. A divisor is {\rm orientable} if it can be obtained this way.

Fix $q\in V(G)$. An orientation $\Ocal$ is {\rm $q$-connected} if any vertex of $G$ can be reached from $q$ via a directed path. A divisor is {\rm $q$-orientable} if it is of the form $D_{\Ocal}$ for some $q$-connected orientation $\Ocal$. (Note that a $q$-orientable divisor is equivalent to the indegree sequence of a root-connected orientation in the sense of \cite{Bernardi_Process}.)
\end{definition}

The following observation reduces many questions related to break divisors to orientable divisors.

\begin{proposition}(\cite[Lemma 3.3]{ABKS}, \cite[Theorem 7.5]{Bac_GO}) \label{Break_equal_qconn}
Fix $q\in V(G)$. The map $D\mapsto D-(q)$ gives a bijection between break divisors and $q$-orientable divisors. In terms of orientations, given a spanning tree $T$ and an orientation of the edges not in $T$, one can get a $q$-connected orientation of $G$ by orienting edges in $T$ away from $q$.
\end{proposition}

Next we mention Gioan's work on the {\em cycle-cocycle reversal system} \cite{Gioan_CCsys} in the language of divisors, as in Backman's paper \cite{Bac_GO}.  Given an orientation $\Ocal$ of a graph $G$, a {\em cycle reversal} reverses all edges of a directed cycle of $\Ocal$ and a {\em cocycle reversal} reverses all edges of a directed cut. These two operations have concise interpretations at the level of divisors:

\begin{proposition}(\cite[Lemma 3.1, Theorem 3.3]{Bac_GO}, \cite[Proposition 4.10, Corollary 4.13]{Gioan_CCsys})
Let $\Ocal,\Ocal'$ be two orientations of a graph $G$. Then $D_{\Ocal}=D_{\Ocal'}$ if and only if $\Ocal'$ can be obtained from $\Ocal$ by a sequence of cycle reversals, and $D_{\Ocal}\sim D_{\Ocal'}$ if and only if $\Ocal'$ can be obtained from $\Ocal$ by a sequence of cycle reversals and cocycle reversals.
\end{proposition}

As a corollary, we can define an equivalence relation $\sim$ on the set $\mathfrak{O}_G$ of orientations of $G$ by setting $\Ocal\sim\Ocal'$ if $\Ocal$ and $\Ocal'$ differ by a sequence of cycle/cocycle reversals; denote by $[\Ocal]$ the equivalence class an orientation $\Ocal$ is in; denote by $[\Ocal]$ the equivalence class an orientation $\Ocal$ is in. The set of such {\em cycle-cocycle reversal classes} is naturally in bijection with $\Pic^{g-1}(G)$:

\begin{proposition}(\cite[Section 5]{Bac_GO})
The map $\tau_G:\mathfrak{O}_G/\!\!\sim\rightarrow \Pic^{g-1}(G)$ given by $[\Ocal]\mapsto [D_{\Ocal}]$ is a well-defined bijection.
\end{proposition}

Backman, generalizing work of Felsner \cite{Felsner}, gives a polynomial time algorithm \cite[Algorithm 7.6]{Bac_GO} which, given a degree $g$ divisor $D$ and a vertex $q$, produces a $q$-connected orientation $\Ocal$ such that $D\sim D_{\Ocal}+(q)$. In particular, if $D$ is itself a break divisor, then $D=D_{\Ocal}+(q)$. It is worth mentioning that the key ingredient underlying Backman's algorithm is the max-flow-min-cut theorem in graph theory, and the algorithm itself uses a maximum flow algorithm.

\subsection{Bernardi's Process and Bijections} \label{sec:Intro_Bernardi}

Bernardi introduced an algorithmic bijection from the set $S(G)$ of spanning trees of a graph $G$ to the set of indegree sequences of $q$-connected orientations \cite{Bernardi_Process}. His bijection uses the notion of {\em ribbon graph} (also known as a combinatorial map or rotation system), which is a finite graph $G$ together with a cyclic ordering of the edges around each vertex. Every embedding of a graph into a closed orientable surface induces a ribbon structure on $G$, where the cyclic ordering of edges around a vertex comes from the orientation of the surface; conversely every ribbon structure on $G$ induces an embedding of $G$ onto a closed orientable surface up to homeomorphism. We refer the reader to \cite[Section 3.2]{Gross_Tucker_TG} for details. We say a ribbon graph is {\em planar} if the surface of the corresponding embedding is the sphere, and whenever we talk about a planar ribbon graph we implicitly assume the graph is a plane graph, i.e. a graph already embedded into $\mathbb{R}^2$ as a geometric object (cf. \cite[Section 4.2]{Diestel_GT}), and the ribbon structure is induced by the counter-clockwise orientation of the plane. The only exception is when we are working with planar duals (see Section \ref{sec:Intro_Duality}), in which case we assume that the ribbon structures of planar duals are induced by the clockwise orientation of the plane.\\

Fix a ribbon structure of a graph $G$ and a {\em starting pair} $(v,f)$, where $f$ is an edge and $v$ is a vertex incident to $f$. For any spanning tree $T$ of $G$, {\em Bernardi process} produces a tour $(v_0,e_0,v_1,\ldots,v_k,e_k)$ of the vertices and edges of $G$. Explicitly, start with $v_0=v,e_0=f$, and in each step, determine $v_{i+1},e_{i+1}$ using $v_i,e_i$ as follows: if $e_i\not\in T$, then set $v_{i+1}=v_i$ and set $e_{i+1}$ to be the next edge of $e_i$ around $v_i$ in the fixed cyclic ordering; otherwise set $v_{i+1}$ to be the other end of $e_i$ and set $e_{i+1}$ to be the next edge of $e_i$ around $v_{i+1}$. The process stops when every edge is traversed exactly twice. Informally, the tour is obtained by walking along edges belonging to $T$ and cutting through edges not belonging to $T$, beginning with $f$ and proceeding according to the ribbon structure, see \cite{Bernardi_Tour}, \cite{Bernardi_Process} for details.\\

We can associate a break divisor $\beta_{(v,f)}(T)$ to such a tour. For each edge $e$ of $G$, $e$ appears in the tour twice as $e_i,e_j$, $i<j$, denote by $\eta(e)$ the vertex $v_i$. Now $\beta_{(v,f)}(T)$ is defined as $\sum_{e\not\in T}(\eta(e))$, i.e. whenever it is the first time we cut through an edge $e\not\in T$, we put a chip at the vertex we cut $e$ from. It is shown by Bernardi (in different terminology) \cite{Bernardi_Process} that for any fixed ribbon structure and starting pair, the Bernardi process induces a bijection between $S(G)$ and the set of break divisors. We call these bijections {\em Bernardi bijections}. Explicit inverses to $\beta_{(v,f)}$ are given by Bernardi \cite{Bernardi_Process} and Baker-Wang \cite{BW_Torsor}. In this paper, we often abuse notation and interpret $\beta_{(v,f)}$ as a map from $S(G)$ to $\Pic^g(G)$, as well as a map from $S(G)$ to the set of partial orientations of $G$, i.e. whenever $e\not\in T$, orient it from $v_j$ to $v_i$. The divisor $\beta_{(v,f)}(T)$ is thus the in-degree sequence of such a partial orientation. As in Proposition \ref{Break_equal_qconn}, one can extend such a partial orientation to a full orientation by orienting edges in $T$ away from some vertex $q$.\\

\begin{figure}[h!]
\begin{center}
    \includegraphics[width=0.45\textwidth]{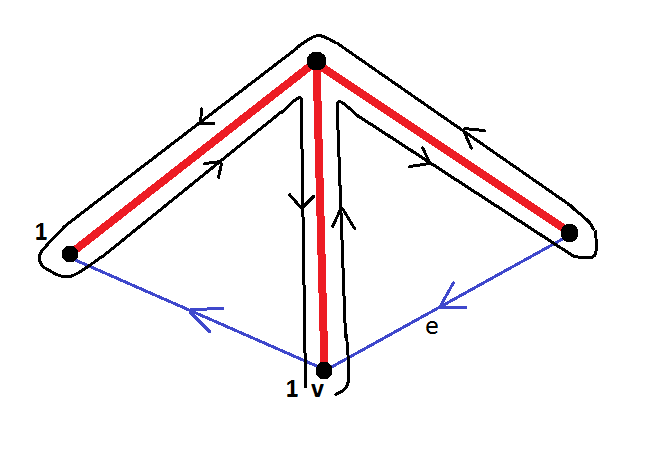}
    \includegraphics[width=0.45\textwidth]{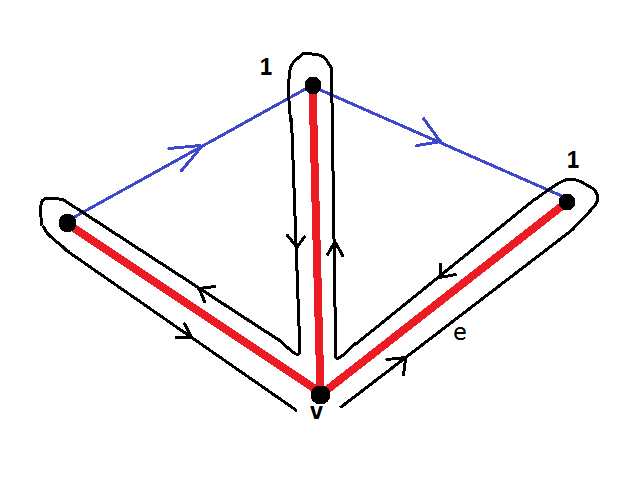}
\end{center}
  \caption{Bernardi processes on two different spanning trees, using the same (planar) ribbon structure and starting pair (\cite[Figure 2]{BW_Torsor}).}
  \label{Bernardi_diag}
\end{figure}

\subsection{Plane Duality and Torsors} \label{sec:Intro_Duality}

Let $G=(V,E)$ be a bridgeless plane graph and let $G^\star=(V^\star,E^\star)$ be a {\em planar dual}, so each face $F$ of $G$ corresponds to a {\em dual vertex} $F^\star$ of $G^\star$, and every edge $e$ of $G$ corresponds to a {\em dual edge} $e^\star$ in which $F_1^\star,F_2^\star$ are adjacent along $e^\star$ if and only if the two faces $F_1,F_2$ are adjacent along $e$ in $G$, see \cite[Chapter 4]{Diestel_GT} for details. Denote by $g^\star:=n-1$ the genus of $G^\star$.\\

\begin{figure}[h!]
\begin{center}
    \includegraphics[width=0.8\textwidth]{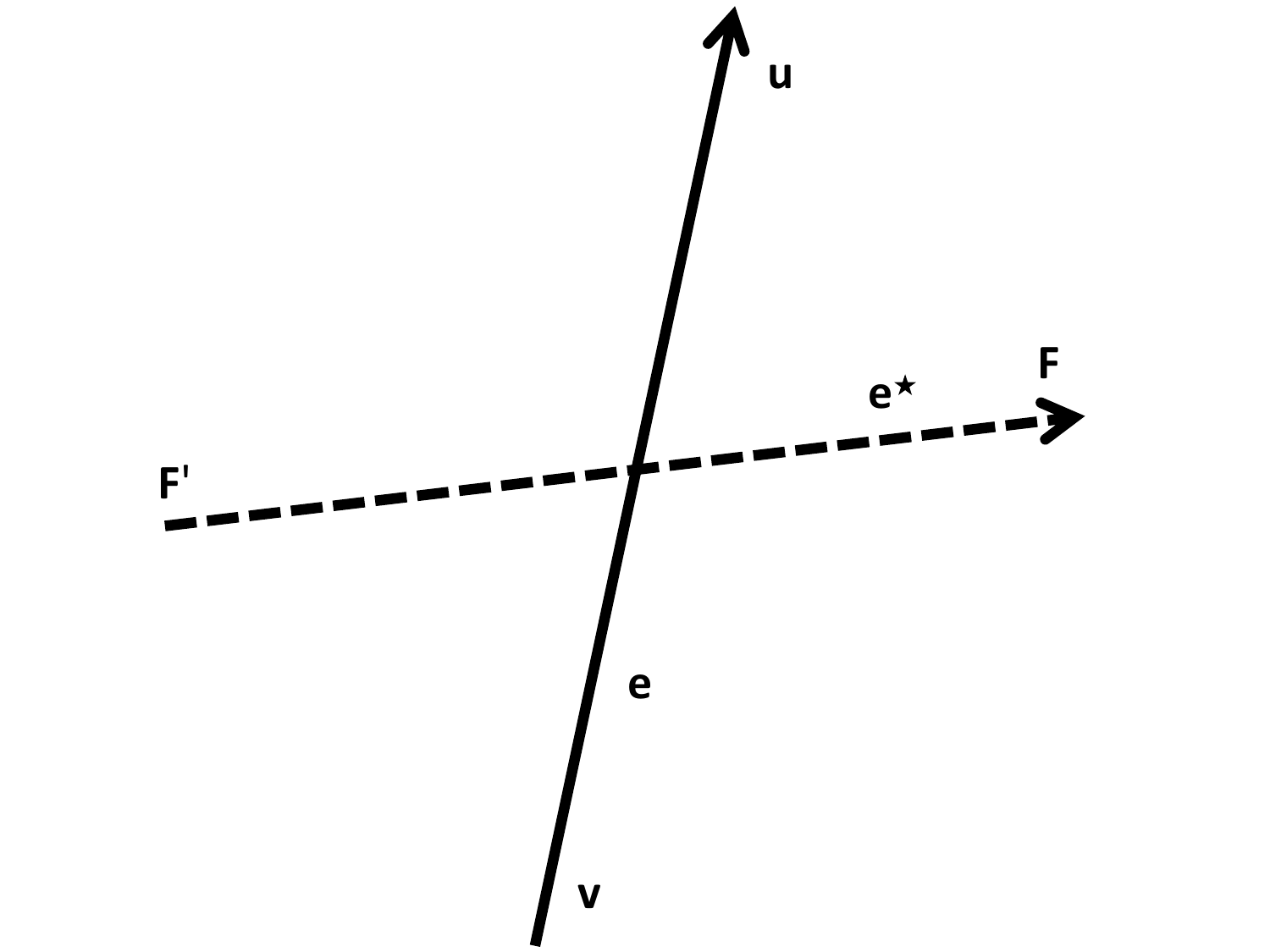}
\end{center}
  \caption{Some conventions on the orientation of plane graphs. Here the orientation of $e$ induces the orientation of the dual edge $e^\star$, and $F$ is said to be the face on the right of $(v,e)$.}
  \label{Plane_diag}
\end{figure}

Every spanning tree $T$ of $G$ corresponds to a {\em dual tree} $T^\star:=\{e^\star:e\not\in T\}$ of $G^\star$ and vice versa, denote the corresponding map from $S(G)$ to $S(G^\star)$ by $\sigma_T$. Every full orientation $\Ocal$ of $G$ corresponds to a {\em dual orientation} $\Ocal^\star$ of $G^\star$, by the convention that locally near the crossing of $e$ and $e^\star$, the orientation of $e^\star$ is obtained from the orientation of $e$ by following the clockwise orientation of the plane, see Figure \ref{Plane_diag} for illustration. Since the directed cycles (resp. cocycles) in $\Ocal$ are the directed cocycles (resp. cycles) in $\Ocal^\star$, the map $\sigma_O:[\Ocal]\mapsto[\Ocal^\star]$ is well-defined and gives a bijection between the cycle-cocycle reversal systems $\mathfrak{O}_G/\!\!\sim,\mathfrak{O}_{G^\star}/\!\!\sim^\star$ of $G$ and $G^\star$, respectively. Composing with the maps $\tau_G,\tau_{G^\star}$ from Section \ref{sec:Intro_GO}, we have a bijective map $\sigma^{g-1}_D:\Pic^{g-1}(G)\rightarrow \Pic^{g^\star-1}(G^\star)$ given by $[D_{\Ocal}]\mapsto [D_{\Ocal^\star}]$. On the other hand, $\Pic^0(G)$ can be identified with $\frac{C_I(G)}{B_I(G)\oplus Z_I(G)}$, where $C_I(G),B_I(G),Z_I(G)$ are the {\em lattice of 1-chains}, the {\em lattice of integer cuts} (cocycles), and the {\em lattice of integer flows} (cycles) of $G$, respectively \cite[Proposition 28.2]{Biggs_APG}. Hence there is a canonical isomorphism $\sigma^0_D:\Pic^0(G)\rightarrow \Pic^0(G^\star)$ using the canonical isomorphism between $B_I(G)$ and $Z_I(G^\star)$, $Z_I(G)$ and $B_I(G^\star)$ \cite[Proposition 8]{BHN_Lattice}, and the identification of $C_I(G)$ and $C_I(G^\star)$.\\

Given a set $X$, a group $G$ and a group action $G\circlearrowright X$, we say $X$ is a {\em $G$-torsor} if the action is simply transitive. For example, $\Pic^{g-1}(G)$ is a $\Pic^0(G)$-torsor via the action $[D]\cdot [E]=[D+E]$ for $[D]\in\Pic^0(G),[E]\in\Pic^{g-1}(G)$. More generally, whenever we have a bijection $\beta:X\rightarrow \Pic^{g-1}(G)$ between some set $X$ and $\Pic^{g-1}(G)$, $X$ is a $\Pic^0(G)$-torsor, where the group action is $[D]\cdot x=\beta^{-1}([D]\cdot \beta(x))$. We will call a torsor induced by a Bernardi bijection {\em Bernardi torsor}. It is routine to check that two bijections $\beta_1,\beta_2:X\rightarrow \Pic^{g-1}(G)$ give {\em isomorphic torsors}, i.e. $\beta_1^{-1}([D]\cdot \beta_1(x))=\beta_2^{-1}([D]\cdot \beta_2(x))$ for all $[D]\in\Pic^0(G),x\in X$, if and only if there exists some translating element $[D_0]\in\Pic^0(G)$ such that $\beta_2(x)=[D_0]+\beta_1(x)$ for all $x\in X$.

\section{Geometric Bijections} \label{sec:GB}

\subsection{Cycle Orientation Maps} \label{sec:COM}

Consider the polyhedral decomposition of $\Pic^g(\Gamma)\cong\Jac(\Gamma)$ introduced in Section \ref{sec:Intro_Decomp}. We can get a bijection between vertices and cells of the decomposition by ``shifting'' \cite[Remark 4.26]{ABKS}. Explicitly, choose some direction ${\bf w}$ in $H_1(G,\mathbb{R})$, superimpose two copies of $\Jac(\Gamma)$, shift one copy along the direction of ${\bf w}$ by an infinitesimally small distance. If ${\bf w}$ is generic (whose definition will be clear below), then each cell in the shifted copy will contain exactly one vertex in the fixed copy (conversely we can say we shift the vertices along the direction $-{\bf w}$). In this manner, we have a bijection $\Psi_{\bf w}$ from cells to vertices, thus a bijection from spanning trees of $G$ to break divisors of $G$.\\

We are going to describe such bijections combinatorially. To do so we first introduce some necessary terminology.

\begin{definition} \label{Def_COC}
A {\rm cycle orientation configuration} is an assignment of an orientation to each cycle of the graph $G$. Given a cycle orientation configuration $\mathfrak{C}$, the {\rm cycle orientation map} corresponding to $\mathfrak{C}$ is a map from the set of spanning trees to the set of break divisors, defined as follows: given a spanning tree $T$ of $G$, orient each edge outside $T$ according to the orientation of its fundamental cycle (specified by $\mathfrak{C}$), then put a chip at the head of each such edge.\\

A vector ${\bf w}$ is {\em generic} if $\la {\bf w}, C\ra\neq 0$ for every cycle $C$ of $G$. Each generic vector ${\bf w}$ induces a cycle orientation configuration $\mathfrak{C}_{\bf w}$ by choosing the orientation of each cycle $C$ in the way such that $\la {\bf w}, C\ra>0$, such cycle orientation configuration and the corresponding cycle orientation map are said to be {\rm geometric}.\\

Given a directed cycle $C$ and an edge $e$, the function $\sign(C,e)$ equals 0 if $e\not\in C$, equals 1 if the orientation of $e$ in $C$ agrees with the reference orientation of $e$, and equals $-1$ otherwise.
\end{definition}

\begin{theorem} \label{main_theorem}
For a generic ${\bf w}$, the geometric bijection $\Psi_{\bf w}$ equals the cycle orientation map corresponding to the cycle orientation configuration $\mathfrak{C}_{\bf w}$. In particular, a geometric cycle orientation map is always bijective.
\end{theorem}

The rest of this section will be devoted to prove Theorem \ref{main_theorem}, as well as describing a complete fan in $\mathbb{R}^g$ that captures all geometric cycle orientation configurations.  We will proceed in three steps: first we work with a single cell as a subset of $H_1(G,\mathbb{R})$, then we put the same cell back to $\Jac(\Gamma)$, and finally we handle all cells in $\Jac(\Gamma)$ together. We use the following convention: whenever we specify a spanning tree $T$ of $G$, we will fix an arbitrary order of the edges outside $T$ and name them as $e_1,\ldots,e_g$, and we write the shifting vector ${\bf w}$ as $\alpha_1\pi(e_1)+\ldots+\alpha_g\pi(e_g)$.\\

For a spanning tree $T$ of $G$, recall that the cell $\Ccal_T$ can be thought as the image of the parallelotope $P_T:=\{\sum_{i=1}^g\theta_i\pi(e_i):0\leq\theta_i\leq 1\}\subset H_1(G,\mathbb{R})$ in $\Jac(\Gamma)=H_1(G,\mathbb{R})/H_1(G,\mathbb{Z})$. We list some elementary and trivial properties of a parallelotope here.

\begin{proposition} \label{trivial_parallel}
Let $P=\{\sum_{i=1}^g \theta_i u_i:0\leq\theta_i\leq 1\}\subset\mathbb{R}^g$ be a full-dimensional parallelotope. For a vector ${\bf w}\in\mathbb{R}^g$, a vertex $v$ of $P$ will be shifted into the interior of $P$ along $-{\bf w}$ if and only if $-{\bf w}$ is in the interior of the {\em tangent cone} $C_v:=\{{\bf u}\in\mathbb{R}^g:\exists\epsilon>0, v+\epsilon{\bf u}\in P_T\}$ of $v$.\\

If the vertex $v$ of $P$ is given by $\sum_{i=1}^g s_i u_i$ (here $s_1,\ldots,s_g\in\{0,1\}$), then $C_v$ is the cone generated by the vectors $(-1)^{s_1} u_1,(-1)^{s_2} u_2,\ldots,(-1)^{s_g} u_g$, i.e. $C_v=\{\sum_{i=1}^g \alpha_i u_i:(-1)^{s_i}\alpha_i\geq 0\}$.\\

If we take the union of tangent cones over all vertices of $P$, then we obtain a fan corresponding to the hyperplane arrangement with hyperplanes $\Lspan\{u_1,\ldots,\widehat{u_i},\ldots,u_g\}$ for $i=1,2,\ldots,g$, here $\widehat{u_i}$ means we omit $u_i$ in the span.
\end{proposition}

We denote the fan in Proposition \ref{trivial_parallel} by $\Fcal_T$ if the parallelotope we are considering is $P_T$. In such case, the hyperplanes in the fan have a simple combinatorial interpretation.

\begin{proposition} \label{HP_interpret}
Let $f_1,\ldots,f_{g-1}$ be distinct edges of $G$ such that $G'=G-\{f_1,\ldots,f_{g-1}\}$ is connected, let $C$ be the unique cycle of $G'$, and let $H=\Lspan\{\pi(f_1),\ldots,\pi(f_{g-1})\}\subset H_1(G,\mathbb{R})$ . Then for an edge $f$ of $G$, $\pi(f)\in H$ if and only if $f\not\in E(C)$. Furthermore, $H$ equals $C^{\perp}$, the orthogonal complement of the vector $C\in H_1(G,\mathbb{R})$.
\end{proposition}

\noindent{\sc Proof:} Given $f\not\in\{f_1,\ldots,f_{g-1}\}$, $\pi(f)\in H$ if and only if $\pi(f_1),\ldots,\pi(f_{g-1})$ and $\pi(f)$ do not span the whole $H_1(G,\mathbb{R})$, if and only if $G-\{f_1,\ldots,f_{g-1},f\}$ is not connected, if and only if $f\not\in E(C)$. For the second part, note that for $i=1,2,\ldots,g-1$, $\la \pi(f_i),C\ra=\la f_i,\pi(C)\ra=\la f_i,C\ra$, which is 0 because $f_i\not\in E(C)$.\hfill $\Box$\\

When $f_1,\ldots,f_{g-1},f$ are the edges outside $T$, the cycle $C$ in Proposition \ref{HP_interpret} is the fundamental cycle of $f$ with respect to $T$. Therefore if we consider Proposition \ref{trivial_parallel} in the context of $P_T$, we see that a (unique) vertex of $P_T$ is shifted into the interior of $P_T$ along $-{\bf w}$ if and only if ${\bf w}$ is not contained in any hyperplane of the form $C^{\perp}$ where $C$ is a fundamental cycle of $T$, if and only if all $\alpha_i$'s are non-zero. Moreover, if all $\alpha_i$'s are indeed non-zero, then the vertex being shifted into $P_T$ is $\sum_{i=1}^g s_i({\bf w})\pi(e_i)$, where $s_i({\bf w})$ is 0 if $\alpha_i<0$ and is 1 otherwise.\\

Now we proceed to the second step. The map from $P_T\subset H_1(G,\mathbb{R})$ to $\Jac(\Gamma)$ might not be injective in general, but recall that the coordinates of a point in $P_T$ with respect to $\pi(e_1),\ldots,\pi(e_g)$ reflect the position of the chips in the edges $e_1,\ldots,e_g$, so non-injectivity occurs precisely when there are more than one way to put chips on those edges to produce the same break divisor. Such ambiguity can not happen in the interior of $P_T$ as all chips are located in the interior of $e_i$'s, hence the map $P_T\rightarrow\Jac(\Gamma)$ is injective when restricted to the interior of $P_T$. Also note that a vertex of $P_T$ will only be mapped to a vertex of $\Ccal_T$, therefore it still makes sense to talk about the fan $\Fcal_T$ of $\Ccal_T$, and the image divisor of $\sum_{i=1}^g s_i({\bf w})\pi(e_i)\in P_T$ discussed above will be the break divisor $\Psi_{\bf w}(T)$.\\

Now we are ready to prove Theorem \ref{main_theorem}.\\

\noindent{\sc Proof of Theorem \ref{main_theorem}:} We need to prove the break divisor $\Psi_{\bf w}(T)$ equals the output of the cycle orientation map on $T$ with respect to $\mathfrak{C}_{\bf w}$. The image of $\sum_{i=1}^g s_i({\bf w}) \pi(e_i)\in P_T$ in $\Jac(\Gamma)$ is the integral break divisor $\sum_{i=1}^g (\eta_{s_i({\bf w})}(e_i))$, where $\eta_{0}(e_i)$ equals the tail of $e_i$ and $\eta_{1}(e_i)$ equals the head of $e_i$. In terms of orientation, for each $i=1,2,\ldots, g$, orient $e_i$ as its reference orientation if $\alpha_i>0$ and as the opposite orientation if $\alpha_i<0$.\\

Let $C_i$ be the fundamental cycle of $e_i$ with respect to $T$, oriented according to $\mathfrak{C}_{\bf w}$. We have $0<\la {\bf w}, C_i\ra=\la\sum_{j=1}^g \alpha_j \pi(e_j), C_i\ra=\sum_{j=1}^g \alpha_j\la \pi(e_j), C_i\ra=\alpha_i\sign(C_i,e_i)$, thus the sign of $\alpha_i$ equals $\sign(C_i,e_i)$. If $\alpha_i>0$, then $e_i$ is oriented according to its reference orientation in $\Psi_{\bf w}(T)$ as discussed in the last paragraph, which is the same as its orientation in $C_i$ as $\sign(C_i,e_i)=1$ here; similarly if $\alpha_i<0$, then $e_i$ is oriented against to its reference orientation in $\Psi_{\bf w}(T)$, which is also the same as its orientation in $C_i$ as $\sign(C_i,e_i)=-1$. \hfill $\Box$\\

In the last step, we take the common refinement of all $\Fcal_T$'s as $T$ varies over all spanning trees and produce the following fan.

\begin{definition} \label{def_GFan}
The {\em geometric bijection fan} of $G$ is the fan constructed by partitioning $H_{1}(G,\mathbb{R})$ using all hyperplanes of the form $C^\perp$, where $C$ varies over all cycles that are the fundamental cycles with respect to some spanning trees.
\end{definition}

A quick observation is that we are actually considering all cycles of $G$ in the definition above.

\begin{proposition} \label{every_cyc_GBCS}
If $C$ is a cycle of a graph $G$, then $C^{\perp}$ occurs as a hyperplane
in the geometric bijection fan of $G$.
\end{proposition}

\noindent{\sc Proof:} Pick any edge $e\in E(C)$. Since $G-e$ is connected and $C-e$ is acyclic, there exists some spanning tree $T$ of $G$ that contains all edges of $C$ except $e$. Therefore $C$ is a fundamental cycle of $T$ and $C^{\perp}$ occurs in the tangent fan of $\Ccal_T$.\hfill $\Box$\\

Therefore in order for $\Psi_{\bf w}$ to be well-defined, ${\bf w}$ must be generic in the sense of Definition \ref{Def_COC}. Finally, there is an obvious bijective correspondence between full-dimensional cones in the geometric bijection fan and geometric cycle orientation configurations: all vectors ${\bf w}$ in the interior of a full-dimensional cone produce the same $\mathfrak{C}_{\bf w}$, while any two vectors from two distinct cones must be on the different sides of some hyperplane $C^\perp$, so the cycle orientation configurations they produced are different.

\subsection{Criteria for Being Geometric}\label{sec:Geom_Criterion}

Not every cycle orientation configuration is geometric. In this section we give two combinatorial criteria for a cycle orientation configuration to be geometric.

\begin{proposition} \label{First_Criterion}
A cycle orientation configuration is geometric if and only if one can assign weights $\alpha_1,\ldots,\alpha_m$ to the edges $e_1,\ldots,e_m$\footnote{More generally, it suffices to assign weights only to edges outside some fixed forest.} such that for each cycle $C$, $\sum\sign(C,e_i)\alpha_i>0$.
\end{proposition}

\noindent{\sc Proof:} This follows from the discussion in Section \ref{sec:COM}. Since
$$\sum_{i=1}^m\sign(C,e_i)\alpha_i=\la \sum_{i=1}^m\alpha_i\pi(e_i),C\ra,$$
a cycle orientation configuration induced by the weights comes from shifting along the vector $ \sum_{i=1}^m\alpha_i\pi(e_i)$. Conversely, if we write the shifting vector as a linear combination of $\pi(e_1),\ldots,\pi(e_m)$, then we can take the coefficients as weights.\hfill $\Box$\\

It is often still difficult to show that a cycle orientation configuration is (not) geometric using Proposition \ref{First_Criterion}. The following alternative criterion is a direct corollary of the Gordan's alternative in linear programming \cite[P. 478]{BLSWZ_Oriented_Matroid}, \cite[Theorem 14]{Border_Alter}, which states that given a matrix $A$, either there exists a vector ${\bf x}$ such that each entry of $A{\bf x}$ is positive, or there exists a non-negative, non-zero solution ${\bf y}$ to the system ${\bf y}^TA=0$, but the two cases can not happen at the same time.

\begin{theorem} \label{generic_LP_criterion}
Let $\mathfrak{C}$ be a cycle orientation configuration and let $C_1,\ldots,C_t$ be the list of simple cycles of $G$, each oriented according to $\mathfrak{C}$. Then $\mathfrak{C}$ is geometric if and only if there exist no non-negative solutions to 
\begin{equation}\label{LP_condition}
\mu_1C_1+\mu_2C_2+\ldots+\mu_tC_t=0
\end{equation}
other than the zero solution\footnote{The condition in Theorem \ref{generic_LP_criterion} is often called the {\em acyclic} condition in the oriented matroid literature, cf. \cite{BLSWZ_Oriented_Matroid}.}. Here we interpret each $C_i$ as a signed sum of edges.
\end{theorem}

\noindent{\sc Proof:} Take $A$ to be the $t\times m$ matrix whose rows are indexed by simple cycles of $G$ and columns are indexed by edges of $G$, and $A_{C,e}$ equals $\sign(C,e)$. Now apply Gordan's alternative. \hfill $\Box$\\

\begin{figure}[h!]
\begin{center}
    \includegraphics[width=0.45\textwidth]{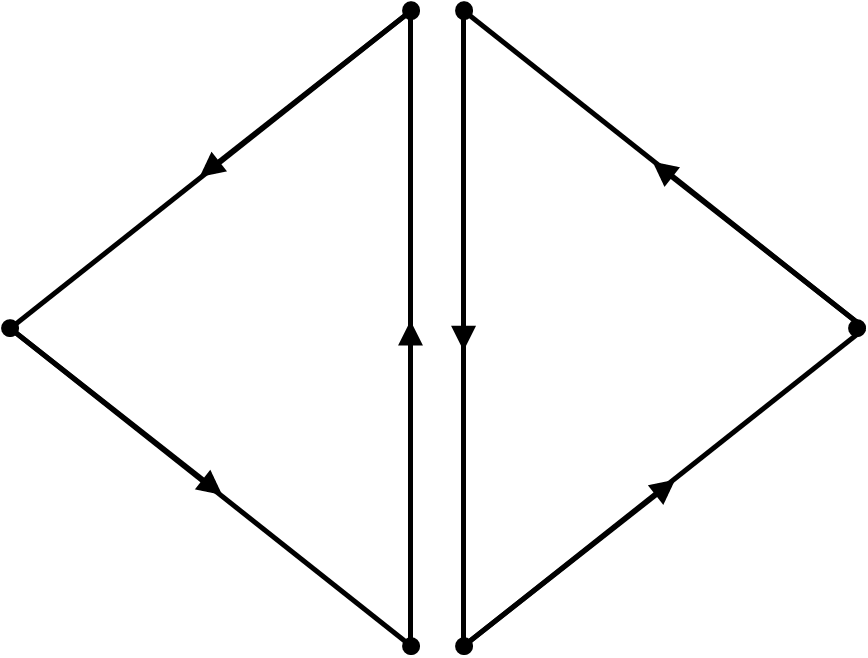}
    \includegraphics[width=0.45\textwidth]{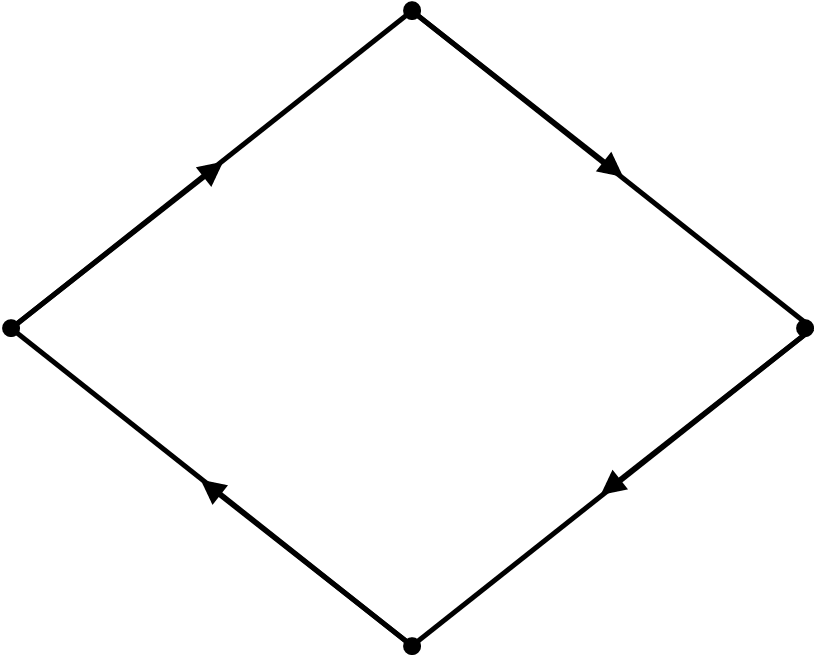}
\end{center}
  \vspace{+5pt}
\begin{center}
    \includegraphics[width=0.45\textwidth]{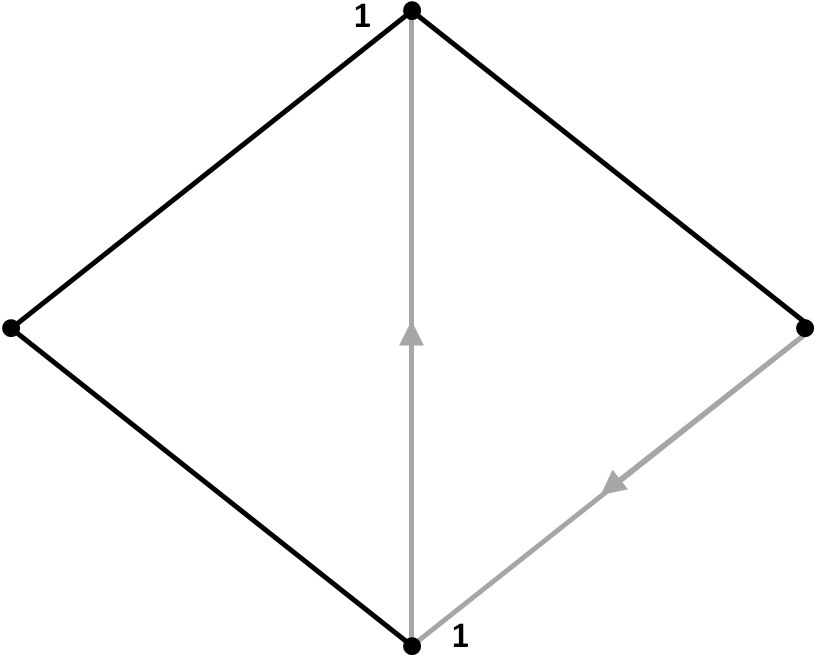}
    \includegraphics[width=0.45\textwidth]{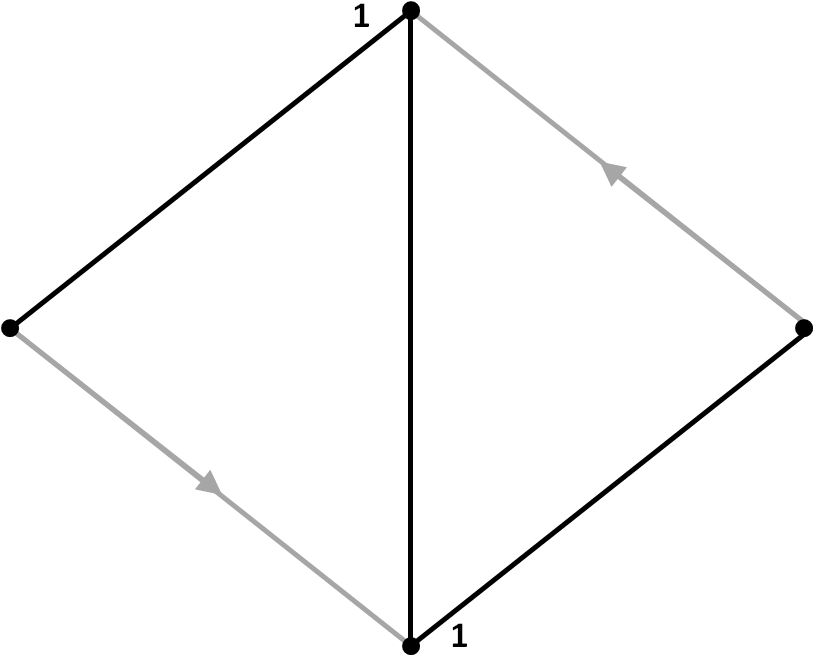}
\end{center}
  \caption{A non-geometric cycle orientation configuration of $\Gamma$ (top), and two distinct spanning trees that map to the same break divisor (bottom).}
  \label{Counter_Example_diag}
\end{figure}

Figure \ref{Counter_Example_diag} shows that not all cycle orientation configurations yield bijective cycle orientation maps, so the acyclic condition (\ref{LP_condition}) is a non-trivial sufficient condition for a cycle orientation map to be bijective.\\

\subsection{The Geometric Bijection Fan}\label{sec:Fan}

We mention two classes of objects that are indexed by the full-dimensional cones of the geometric bijection fan (equivalently, by geometric cycle orientation configurations), which might be of independent interests.

\begin{definition} \label{Def_Compatible}
Let $\mathfrak{C}$ be a cycle orientation configuration of a graph $G$. An orientation $\Ocal$ is {\rm $\mathfrak{C}$-compatible} if every directed cycle of $\Ocal$ is oriented as in $\mathfrak{C}$.
\end{definition}

\begin{proposition} \label{Def_FO}
For each geometric cycle orientation configuration $\mathfrak{C}$ of $G$, the set of $\mathfrak{C}$-compatible orientations form a system of representatives for the cycle reversal system of $G$.
\end{proposition}

\noindent{\sc Proof:} For uniqueness, note that any two distinct orientations in the same cycle reversal class differ by a disjoint union of directed cycles, so they cannot both be $\mathfrak{C}$-compatible. For existence, start with an orientation $\Ocal$ and keep reversing some directed cycle not oriented as in $\mathfrak{C}$, if any. This process will eventually stop: suppose not, since the number of orientations is finite, WLOG the orientation returns to $\Ocal$ after flipping some directed cycles $C_1,\ldots,C_k$ in that order (the cycles might not be distinct). Since $C_1+\ldots+C_k=0$ and each $C_i$ is of the opposite orientation as in $\mathfrak{C}$, we have a contradiction with (\ref{LP_condition}). \hfill $\Box$\\

In particular, each full-dimensional cone of the geometric bijection fan induces a nice system of representatives for the cycle reversal system of $G$.\\

In \cite{Bac_PO}, \cite{Bernardi_Process}, the authors considered the notion of {\em cycle minimal orientations}: given an ordering $e_1<\ldots<e_m$ of the edges of $G$, together with a reference orientation, we say an orientation $\Ocal$ is {\em cycle minimal} if every directed cycle is oriented according to the reference orientation of its minimum edge. Consider the cycle orientation configuration $\mathfrak{C}_<$ that orients each cycle according to its minimum edge. $\mathfrak{C}_<$ can be proven to be geometric (cf. Theorem \ref{Def_EOM} below), hence the notion of $\mathfrak{C}$-compatibility generalizes the notion of cycle minimality. The significance of such an observation is that the purely combinatorial notion of cycle minimal orientations has a natural generalization in terms of polyhedral geometry (an idea first briefly introduced in \cite{GZ_GO}), this suggests the possibility of generalizing other classical notions in combinatorics related to (edge) orderings to our setup.\\

\noindent{\bf Remark.}  An alternative proof of Proposition \ref{Def_FO} can be obtained from the discussion below, using the unique remainder property of division by a Gr\"{o}bner basis, cf. \cite[Section 4]{Bac_PO}.\\

The second class of objects we mention is the set of monomial initial ideals of the Lawrence ideal $\Tcal_G$ associated to the cycle lattice $H_1(G,\mathbb{Z})$ of $G$. This ideal is also the ideal of the classical quasisymmetry model in statistics \cite{KMS_QS}. Formally, $\Tcal_G\subset K[x_1,\ldots,x_m,y_1,\ldots,y_m]$ is generated by binomials $\phi({\bf v})-\phi(-{\bf v})$ for ${\bf v}\in H_1(G,\mathbb{Z})$, where $\phi(a_1{\bf e}_1+\ldots+a_m{\bf e}_m)=\prod_{a_i>0}x_i^{a_i}\prod_{a_i<0}y_i^{-a_i}$.

\begin{proposition} \label{Law_Ideal}
The collection of monomial initial ideals of $\Tcal_G$ are exactly ideals of the form $\langle\phi(C):C\in\mathfrak{C}\rangle$ as $\mathfrak{C}$ ranges over all geometric cycle orientation configurations of $G$.
\end{proposition}

\noindent{\sc Proof:} It is known that $\{\phi(C)-\phi(-C):C{\rm\ is\ a\ directed\ cycle}\}$ is a (minimal) universal Gr\"{o}bner basis of $\Tcal_G$ \cite[Lemma 3.1]{KMS_QS}, \cite[Proposition 7.8]{MS_Ideal}. Therefore to specify a monomial initial ideal of $\Tcal_G$, it suffices to specify the initial term of each basis element (with respect to the corresponding monomial term order), which is equivalent to choosing an orientation for each cycle of $G$; conversely, every monomial initial ideal specifies a way to pick an orientation for each cycle of $G$ by considering its minimal set of generators. Using a standard fact in Gr\"{o}bner theory \cite[Theorem 1.11]{Stur_GBCP}, each monomial initial ideal $I$ of $\Tcal_G$ is equal to ${\rm in}_{\bf\hat{w}}(\Tcal_G)$ for some weight vector ${\bf\hat{w}}\in\mathbb{R}^{2m}$, but the monomial term order $<_{\hat{w}}$ picks the same set of directed cycles as the weight vector ${\bf w}\in\mathbb{R}^m$ (here the $i$-th coordinate of ${\bf w}$ is equal to ${\bf\hat{w}}_{x_i}-{\bf\hat{w}}_{y_i}$) does in the sense of Section \ref{sec:COM}, hence must be geometric. Conversely, every geometric cycle orientation configuration $\mathfrak{C}$ is induced by some weight vector ${\bf w}$, and the weight vector $({\bf w}\ {\bf 0})\in\mathbb{R}^{2m}$ induces the monomial initial ideal $\langle\phi(C):C\in\mathfrak{C}\rangle$. \hfill $\Box$\\

Geometrically, this is to say that the geometric bijection fan is the quotient of the Gr\"{o}bner fan of $\Tcal_G$ modulo a $(2m-g)$-dimensional lineality space.\\

\noindent{\bf Remark.} The notions considered in this section have their dual versions in which we replace ``cycles'' with ``{\em cocycles}'' (also known as minimal cuts or bonds in other literature). Namely, we have $\mathfrak{C}^*$-compatible orientations with respect to a cocycle orientation configuration $\mathfrak{C}^*$ and the Lawrence ideal associated to the cocycle lattice, and our results have corresponding dual counterparts. In particular, given a geometric cycle orientation configuration $\mathfrak{C}$ and a geometric cocycle orientation configuration $\mathfrak{C}^*$, we say an orientation is {\em $(\mathfrak{C},\mathfrak{C}^*)$-compatible} if each directed (co)cycle is oriented as in $\mathfrak{C}$ (resp. $\mathfrak{C}^*$), this generalizes Backman's notion of cycle-cocycle minimal orientations.

\section{A New Combinatorial Bijection} \label{sec:Bijection}

\subsection{Edge Ordering Maps} \label{sec:EOM}

By finding vectors that satisfy some stronger ``genericity'' conditions, we have a surprisingly simple family of geometric cycle orientation configurations/geometric bijections. But to the best of our knowledge, even such special cases are new in the literature.

\begin{theorem} \label{Def_EOM}
Order the edges of $G$ as
$e_{1}<e_{2}<\ldots <e_{m}$ and pick an arbitrary reference orientation for them.\footnote{As in Proposition \ref{First_Criterion}, it suffices to work with edges 
outside some fixed forest.} For each cycle $C$ of the graph, orient $C$ according to the smallest edge contained in $C$. Then the resulting cycle orientation configuration is geometric.
\end{theorem}

\noindent{\sc Proof:} Set the weight $\alpha_i:=\frac{1}{2^i}$ for each edge $e_i$, $i=1,2,\ldots, m$. Any (non-empty) signed subset sum of $\alpha_i$'s is non-zero, so {\em a priori} the signed sum over any cycle is non-zero. Moreover, the sign of a sum depends solely on the sign of the largest $\alpha_i$ appearing in the sum, so it suffices to orient $C$ according to the smallest edge to guarantee that the signed sum is positive.\hfill $\Box$\\

We call the geometric bijections arising from edge orderings as in Theorem \ref{Def_EOM} as {\em edge ordering maps/bijections}. As explained in the proof of Proposition \ref{First_Criterion}, each edge ordering map corresponds to the cone in the geometric bijection fan containing the vector $\sum_{i=1}^m\frac{1}{2^i}\pi(e_i)$ in its interior.\\

Assuming the initial data from Theorem \ref{Def_EOM}, a pseudocode for edge ordering map is given below.

\begin{algorithm}\caption{Edge Ordering Map}
\KwIn{A spanning tree $T$ of $G$.}
\KwOut{A break divisor $D\in\Div(G)$.}
\BlankLine
Set $D:=0$.
\BlankLine
\For{$f\not\in T$}{
	$C:={\rm Unique\ cycle\ contained\ in\ }T+f$\;
	$i:={\rm index\ of\ the\ smallest\ edge\ in\ }C$\;
	Orient $f$ to have the same orientation as $e_i$ in $C$\;
	$D:=D+{\rm Head}(f)$\;
}
\BlankLine
Output $D$.
\label{EdgeOrderingMap}
\end{algorithm}

\subsection{An Inverse Algorithm} \label{sec:Inverse}

In this section, we give a combinatorial inverse algorithm for Algorithm \ref{EdgeOrderingMap}, thereby providing a combinatorial proof that the map given by Algorithm \ref{EdgeOrderingMap} is indeed a bijection.\\

\noindent{\bf Remark.} As the example in the Appendix shows, the notion of geometric cycle orientation configurations is strictly more general than those configurations coming from edge orderings, hence the combinatorial proof presented here is not amenable to the general case.\\

We will first describe the key subroutine {\em Inverse} in Algorithm \ref{Inverse_Routine}, which works at the level of orientations; then we will give the main algorithm in Algorithm \ref{Inverse_EOM}. Here the subroutine {\em DivisorToOrientation} is the algorithm by Backman \cite[Algorithm 7.6]{Bac_GO} which, given a break divisor $D$ on a graph $G$ and a vertex $q$, outputs a $q$-connected orientation $\Ocal$ such that $D=D_{\Ocal}+(q)$.

\begin{algorithm}
\caption{The subroutine Inverse}
\KwIn{A connected graph $H$, a vertex $q\in V(H)$ and a $q$-connected orientation $\Ocal$.}
\KwOut{A spanning tree $T$ of H.}
\BlankLine
\eIf{$H$ is a tree}{
$T:=H$;}{
$i:={\rm index\ of\ the\ smallest\ edge\ in\ }H$;\\
$P:=$ A directed path from $q$ to $v_i$ in $\Ocal$;\\
Reverse the edges of $P$ in $\Ocal$ to obtain $\hat{\Ocal}$;\\
$U:={\rm Vectices\ }u_i{\rm\ can\ reach\ in\ }\hat{\Ocal}$;\\
	\eIf{$v_i\in U$}{
	  \If{$e_i$ goes from $v_i$ to $u_i$ in $\hat{\Ocal}$}{
	    $Q:=$ A directed path from $u_i$ to $v_i$ in $\hat{\Ocal}$;\\
	    Reverse the edges of the directed cycle $\{e_i\}\cup Q$ in $\hat{\Ocal}$;
	  }
	$T:={\rm Inverse}(H-e_i,v_i,\hat{\Ocal}-e_i)$;\\}{
	$T:=\{e_i\}\cup{\rm Inverse}(H[U],u_i,\hat{\Ocal}|_{H[U]})\cup{\rm Inverse}(H[U^{c}],v_i,\hat{\Ocal}|_{H[U^c]})$;
	}
}
\BlankLine
Output $T$.
\label{Inverse_Routine}
\end{algorithm}

\begin{algorithm}
\caption{Inverse to the Edge Order Map}
\KwIn{A connected graph $G$ and a break divisor $D\in\Div(G)$.}
\KwOut{A spanning tree $T$ of G.}
\BlankLine
$\Ocal:={\rm DivisorToOrientation}(G,D,v_1)$;\\
$T:={\rm Inverse}(G,v_1,\Ocal)$;
\BlankLine
Output $T$.
\label{Inverse_EOM}
\end{algorithm}

\begin{theorem}
Algorithm \ref{Inverse_EOM} always terminates and it is the inverse of Algorithm \ref{EdgeOrderingMap}.
\end{theorem}

\noindent{\sc Proof:} We shall first prove every recursive step in Inverse is valid, then show by induction on the number of edges that Inverse is the inverse of Algorithm \ref{EdgeOrderingMap}, namely if $T={\rm Inverse}(H,q,\Ocal)$, then the break divisor associated to $T$ via Algorithm \ref{EdgeOrderingMap} is $D_{\Ocal}+(q)$.\\

The case of $H$ being a tree is obviously correct. For the non-trivial cases, first notice that $D_{\hat{\Ocal}}+(v_i)=D_{\Ocal}+(q)$, so they correspond to the same break divisor. Furthermore, $\hat{\Ocal}$ is $v_i$-connected: for any vertex $u$ in $H$, pick a directed path $R$ from $q$ to $u$ in $\Ocal$, let $w$ be the vertex in the intersection of $V(P)$ and $V(R)$ that is closest to $v_i$ on $P$ ($w$ could be $q$ itself), then concatenating the portion of $P$ from $v_i$ to $w$ and the portion of $R$ from $w$ to $u$ gives a directed path from $v_i$ to $u$ in $\hat{\Ocal}$.\\

In the case of $v_i\in U$, the algorithm performs a cycle reversal if necessary to guarantee $e_i$ goes from $u_i$ to $v_i$ in $\hat{\Ocal}$, so $\hat{\Ocal}-e_i$ is a $v_i$-connected orientation of $H-e_i$ as $v_i$ never needed to use $e_i$ to reach any other vertices in the new $\hat{\Ocal}$, Hence the recursive call ${\rm Inverse}(H-e_i,v_i,\hat{\Ocal}-e_i)$ is valid. Letting $T={\rm Inverse}(H-e_i,v_i,\hat{\Ocal}-e_i)$, by induction the break divisor obtained from $T$ in $H-e_i$ using Algorithm \ref{EdgeOrderingMap} is $D_{\hat{\Ocal}-e_i}+(v_i)$. Consider the fundamental cycle $C$ of $e_i$ in $T+e_i\subset H$. Since $e_i$ is the smallest edge in $C$, $e_i$ will be oriented as its own reference orientation $\overrightarrow{u_iv_i}$ in Algorithm \ref{EdgeOrderingMap}, hence the break divisor obtained from $T$ in $H$ using Algorithm \ref{EdgeOrderingMap} is indeed $D_{\hat{\Ocal}-e_i}+2(v_i)=D_{\hat{\Ocal}}+(v_i)$.\\

In the last case with $v_i\not\in U$, every edge between $U$ and $U^c$ goes from $U^c$ to $U$ in $\hat{\Ocal}$. In particular $D_{\hat{\Ocal}}=D_{\hat{\Ocal}|_{H[U]}}+D_{\hat{\Ocal}|_{H[U^c]}}+\sum_{e=uu',u\in U,u'\in U^c}(u)$, here we abuse notations and consider $D_{\hat{\Ocal}|_{H[U]}}$ and $D_{\hat{\Ocal}|_{H[U^c]}}$ as divisors on $H$ in the obvious way. On one side, $\hat{\Ocal}$ restricted to $H[U^c]$ is $v_i$-connected as $v_i$ could never use edges between $U$ and $U^c$ to access vertices in $U^c$, hence ${\rm Inverse}(H[U^{c}],v_i,\hat{\Ocal}|_{H[U^c]})$ is a valid call; on the other side, $\hat{\Ocal}$ restricted to $H[U]$ is $u_i$-connected by the construction of $U$, hence ${\rm Inverse}(H[U],u_i,\hat{\Ocal}|_{H[U]})$ is also a valid call. Suppose $T$ is the outputted tree here, we consider the break divisor associated to $T$ via Algorithm \ref{EdgeOrderingMap}. By induction hypothesis, the contribution of those non-tree edges in $H[U]$ and $H[U^c]$ is equal to $D_{\hat{\Ocal}|_{H[U]}}+(u_i)+D_{\hat{\Ocal}|_{H[U^c]}}+(v_i)$. For every edge $f\neq e_i$ between $U$ and $U^c$, the fundamental cycle of $T+f$ contains $e_i$, so $f$ will be oriented from $U^c$ to $U$, thus these edges contribute $[\sum_{e=uu',u\in U,u'\in U^c}(u)]-(u_i)$ to the final divisor. Summing these contributions of non-tree edges gives $D_{\hat{\Ocal}}+(v_i)$ as claimed.\hfill $\Box$\\

Now we analyze the complexity of Algorithm \ref{Inverse_EOM}. The complexity of the subroutine DivisorToOrientation is essentially the complexity of a maximum flow algorithm (here any exact algorithm for unit capacity directed graphs suffices, say the $\tilde{O}(m^{10/7})$ algorithm by Madry \cite{Madry_MF}). Each instance of the first case of Inverse takes $O(|V(H)|)$ time, but any two $H$'s considered in the computational process of Algorithm \ref{Inverse_EOM} are disjoint, so the first case takes $O(n)$ time in total. Lastly, for each non-trivial call of Inverse, a different smallest edge $e_i$ is being considered, so there can be $O(m)$ such calls, and it is easy to see such a call can be handled in $O(m)$ time ($P$, $U$, and $Q$ can all be found by BFS, and other maintenance/modification operations also take $O(m)$ time). Therefore the total time complexity is $O(m^2)$.\\

As an application of our algorithm, one can use it to design a polynomial time sampling algorithm for sampling random spanning trees uniformly exactly, using an information theoretical minimum number of random bits. The idea is that one can generate a random element of $\Pic^g(G)$ in polynomial time using only $\lceil\log_2 |S(G)|\rceil$ random bits (cf. \cite[Algorithm 5]{BS_Chip}), and using Backman's algorithm \cite{Bac_GO} one can find the break divisor in the random divisor class of $\Pic^g(G)$. Finally, using our inverse algorithm, one can output a random spanning tree of $G$. We refer the reader to the paper by Baker and Shokrieh \cite{BS_Chip} for further discussion of this paradigm, and \cite{GRV_Random} for applications of random spanning trees sampling algorithms.

\section{Comparison of Geometric Bijections and Bernardi Bijections} \label{sec:Bernardi}

\subsection{Every Planar Bernardi Bijection is Geometric} \label{sec:Plane_Bernardi}

In this section we prove that every Bernardi bijection on a plane graph is a geometric bijection, which was a conjecture by Baker \cite[Remark 5.2]{BW_Torsor}. We adopt the convention that the ribbon structures on plane graphs are induced by the counter-clockwise orientation of the plane.

\begin{proposition}\label{Plane_BP_as_OM}
Let $G$ be a planar ribbon graph. Then every Bernardi process on $G$ is a cycle orientation map. More explicitly, if the Bernardi process starts with the pair $(x,e=xy)$, then the cycles are oriented as follows: if $e$ is outside $C$, orient $C$ clockwise; if $e$ is inside $C$, orient $C$ counter-clockwise; if $e$ is on $C$, orient $C$ as the opposite of $\overrightarrow{xy}$.
\end{proposition}

\begin{figure}[h!]
\begin{center}
    \includegraphics[width=0.8\textwidth]{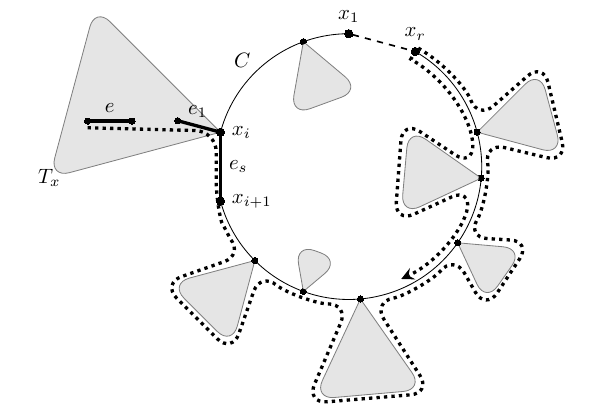}
\end{center}
  \caption{Illustration of Proposition \ref{Plane_BP_as_OM}.}
  \label{Bernardi_OM_diag}
\end{figure}

\noindent{\sc Proof:} Let $C=x_1x_2\ldots x_r$ be a cycle with vertices indexed in counter-clockwise order, and let $T$ be a spanning tree for which $C$ is a fundamental cycle. WLOG $T$ is missing the edge $\hat{e}:=x_rx_1$ from $C$. Suppose $e$ is outside $C$. Then the Bernardi tour starting from $e$, when restricted to $C$, traverses the ``outside'' of $C$ in a counter-clockwise manner before going to the ``inside'' of $C$ by cutting through $\hat{e}$. Hence the tour will put a chip at $x_r$, which corresponds to orienting $C$ clockwise. See Figure \ref{Bernardi_OM_diag} for illustration. The analysis of the remaining cases are similar. \hfill $\Box$\\

One immediate observation is that a more ``correct'' way to index Bernardi processes/bijections on a plane graph is by faces rather than by edges.

\begin{corollary}\label{Face_BP_Bijection}
Let $G$ be a planar ribbon graph and let $(v,e=uv)$ be a starting pair. Let $F$ be the face to the right of $(v,e)$ (cf. Figure \ref{Plane_diag}). Then the Bernardi bijection $\beta_{(v,e)}$ is equal to the cycle orientation map corresponding to the following cycle orientation configuration: a cycle is oriented into counter-clockwise if and only if $F$ is in the interior of the cycle. Conversely every face corresponds to some Bernardi bijection in such manner.
\end{corollary}

For each face $F$ of a plane graph $G$, denote by $\beta_F$ the cycle orientation map coming from the cycle orientation configuration described in Corollary \ref{Face_BP_Bijection}, so every Bernardi bijection on $G$ is a $\beta_F$ for some face $F$ and vice versa.

\begin{algorithm}\caption{From Combinatorial Map to Edge Ordering}
\KwIn{A plane graph $G$.}
\KwOut{An orientation and ordering $e_1<\ldots<e_m$ of edges of $G$.}
\BlankLine
Set $G_1:=G, i:=1$.\\
\BlankLine
\While{There are edges in $G_i$}{
	Pick an edge $e\in G_i$ incident to the unbounded face;\\
	$e_i:=e$;\\
	\eIf{$e_i$ is a cut edge}{
	Orient $e_i$ arbitrarily;}{
	Pick any cycle $C_i\subset G_i$ containing $e_i$;\\
	Orient $e_i$ according to its orientation on a clockwise-oriented $C_i$;\\}
	$G_{i+1}:=G_{i}-e_i$;\\
	$i:=i+1$;
}
\BlankLine
Output $e_1<\ldots<e_m$ and their orientations.
\label{MapToOrder}
\end{algorithm}

Now we prove the main result in this section.

\begin{theorem} \label{Plane_BP_are_geom}
Let $G$ be a planar ribbon graph. Then every Bernardi bijection on $G$ is an edge ordering map, hence a geometric cycle orientation map.\\

Moreover, if $\beta_{F},\beta_{F'}$ are two Bernardi bijections corresponding to adjacent faces $F,F'$ of $G$, then the data describing the two edge ordering maps can be chosen so that they are the same except the orientation of the first edge $e_1$, and $e_1$ can further be chosen to be a common edge of $F,F'$.
\end{theorem}

\noindent{\sc Proof:} Consider Algorithm \ref{MapToOrder} (the algorithm always terminates, as in any step one can find a suitable $e_i$ and orient it). Let $C$ be an arbitrary cycle in $G$, and let $e_i$ be the smallest edge $C$ contains. Then $C\subset G_i$ and $e_i$ is not a cut edge of $G_i$. By construction, $e_i$ is incident to the unbounded face of $G_i$ and inherits its orientation from a clockwise-oriented $C_i$, so orienting $C$ according to $e_i$ will make $C$ clockwise as well. Therefore the output of Algorithm \ref{MapToOrder} will orient every cycle clockwise. This is the cycle orientation configuration for $\beta_F$, where $F$ is the unbounded face. Thus $\beta_F$ is an edge ordering map, and we know from Theorem \ref{Def_EOM} that every edge ordering map is geometric. Our proof (resp. algorithm) can be easily modified for the general case.\\

For the second half of the theorem, note that Algorithm \ref{MapToOrder} can choose the same edge as $e_1$ for both maps, namely a common edge of $F,F'$, but with the opposite orientation. Once $e_1$ is removed, the two cases are the same as $F$ and $F'$ become the same face of $G-e_1$. \hfill $\Box$\\

\noindent{\bf Remark.} A more concise proof to show that a planar Bernardi bijection is geometric can be obtained by checking the condition (1) in Theorem \ref{generic_LP_criterion} directly.

\subsection{The Bernardi Torsors of Plane Graphs and Plane Duality} \label{sec:Plane_torsor}

Baker and Wang proved several propositions about Bernardi torsors of plane graphs in their paper \cite{BW_Torsor}. (They also use the convention that the ribbon structures on planar duals are induced by the clockwise orientation of the plane.) We recall that $\sigma_T,\sigma^0_D,\sigma^{g-1}_D,\sigma_O$ are the canonical duality maps between $S(G)$ and $S(G^\star)$, $\Pic^0(G)$ and $\Pic^0(G^\star)$, $\Pic^{g-1}(G)$ and $\Pic^{g^\star-1}(G^\star)$, and $\mathfrak{O}_G/\!\!\sim$ and $\mathfrak{O}_{G^\star}/\!\!\sim^\star$, respectively; and $\tau_G$ (resp. $\tau_{G^\star}$) is the bijective map between $\mathfrak{O}_G/\!\!\sim$ and $\Pic^{g-1}(G)$ (resp. $\mathfrak{O}_{G^\star}/\!\!\sim^\star$ and $\Pic^{g^\star-1}(G^\star)$).

\begin{theorem} (\cite[Theorem 5.1]{BW_Torsor}) \label{BW_Isom}
Let $G$ be a planar ribbon graph. Then all Bernardi torsors are isomorphic.
\end{theorem}

\begin{theorem} (\cite[Theorem 6.1]{BW_Torsor}) \label{BW_Dual}
Let $G$ be a bridgeless planar ribbon graph, and let $\beta,\beta^\star$ be some Bernardi bijections on $G$ and $G^\star$, respectively.  Then the following diagram is commutative:
\[
\begin{CD}
\Pic^0(G) \times S(G) @>{\beta}>> S(G) \\
@VV{\sigma^0_D \times \sigma_T}V	@VV{\sigma_T}V \\
\Pic^0(G^\star) \times S(G^\star) @>{\beta^\star}>> S(G^\star) \\
\end{CD}
\]
\end{theorem}

\begin{theorem} (\cite[Theorem 7.1]{BW_Torsor}) \label{BW_Rotor}
Let $G$ be a planar ribbon graph.  Then the Bernardi and rotor-routing processes define the same $\Pic^0(G)$-torsor structure on $S(G)$.
\end{theorem}

We give new proofs to the first two propositions using our language of cycle orientation maps and edge ordering maps. The new proofs are more concise and produce more precise versions of statements. Combined with Theorem \ref{BW_Rotor}, our arguments give a new proof that the torsors induced by rotor-routing on a plane graph are all isomorphic and respect plane duality.

\begin{theorem}\label{EOM_torsor}
Let $\beta,\beta '$ be edge ordering maps with data $e_1<e_2<\ldots<e_m$ and $\overline{e}_1<e_2<\ldots<e_m$, respectively, where $\overline{e}_1$ denotes a reversed $e_1$. Then the two bijections give isomorphic torsors, with translating element $[\partial(e_1)]\in \Pic^0(G)$.
\end{theorem}

\noindent{\sc Proof:} Given a spanning tree $T$, we compute $\beta(T)-\beta'(T)$. Say $e_1=\overrightarrow{vu}$. If $e_1\not\in T$, then the only edge that changes orientation is $e_1$ so the difference is $(u)-(v)=\partial(e_1)$. Otherwise $e_1\in T$ and let $K=[V',U']$ be the fundamental cut of $e_1$ with $v\in V',u\in U'$. An edge (not in $T$) changes orientation if and only if its fundamental cycle contains $e_1$, if and only if it is in $K$. Edges in $K\setminus\{e_1\}$ orient from $U'$ to $V'$ in the first edge ordering map and from $V'$ to $U'$ in the second map, so the difference is $\sum_{f=u'v'\in K\setminus\{e_1\},u'\in U',v'\in V'}(v')-(u')=(u)-(v)-\Delta(\rchi_{U'})$, which is equal to $(u)-(v)$ in $\Pic^0(G)$.\hfill $\Box$

\begin{corollary} (More precise version of Theorem \ref{BW_Isom})
Let $G$ be a planar ribbon graph. Then all Bernardi torsors are isomorphic. More precisely, suppose two bijections correspond to the two faces $F,F'$ of $G$, and let $F^\star=F_0^\star-f_1^\star-F_1^\star-f_2^\star-\ldots-F_l^\star=F'^\star$ be a trail in $G^\star$. Then the translating element between the two induced torsors is $[\partial(f_1+f_2+\ldots+f_l)]\in Pic^0(G)$ (with a suitable orientation of the $f_i$'s).
\end{corollary}

\noindent{\sc Proof:} This follows easily by repeatedly applying Theorem \ref{EOM_torsor} to Theorem \ref{Plane_BP_are_geom}. \hfill $\Box$\\

Now we give an alternative proof of Theorem \ref{BW_Dual}, the idea is to prove the commutativity of two finer diagrams separately. In particular, our proof produces a stronger assertion than the proof in \cite{BW_Torsor}.

\begin{proposition} \label{Pic_CD}
Let $G$ be a bridgeless planar ribbon graph and let $G^\star$ be a planar dual of $G$. Then the following diagram commutes, here the horizontal arrows correspond to the addition map:
\[
\begin{CD}
\Pic^0(G) \times\Pic^{g-1}(G) @>>> \Pic^{g-1}(G)\\
@VV{\sigma^0_D\times\sigma^{g-1}_D}V	@VV{\sigma^{g-1}_D}V\\
\Pic^0(G^\star) \times\Pic^{g^\star-1}(G^\star)@>>> \Pic^{g^\star-1}(G^\star)\\
\end{CD}
\]
\end{proposition}

\noindent{\sc Proof:} Since the graph is connected, $\Pic^0(G)$ is generated by elements of the form $[(v)-(u)]$, where $u,v$ are adjacent vertices. Hence by linearity it suffices to prove
$$\sigma^0_D([(v)-(u)])+\sigma^{g-1}_D([D])=\sigma^{g-1}_D([(v)-(u)]+[D])$$
for these $[(v)-(u)]$'s. Fix two such adjacent vertices $u,v$, say they are both incident to the edge $e$. For a divisor class $[D]\in\Pic^{g-1}(G)$, we can interpret the addition $[(v)-(u)]+[D]$ as follows: pick an orientation $\Ocal$ such that $D\sim D_{\Ocal}$ and that $e$ is oriented from $v$ to $u$ (the latter is always possible by reversing a directed cycle/cut whose $e$ is in if necessary), reverse $e$ in $\Ocal$ to obtain $\Ocal'$, then $[(v)-(u)]+[D_{\Ocal}]=[D_{\Ocal'}]$. Using the naming convention in Figure \ref{Plane_diag}, the dual element of $[(v)-(u)]$ is $[(F'^\star)-(F^\star)]\in\Pic^0(G^\star)$, and the dual element of $[D]=[D_{\Ocal}]$ is $[D_{\Ocal^\star}]\in\Pic^{g^\star-1}(G^\star)$, note that $e^\star$ is oriented from $F'^\star$ to $F^\star$ in $\Ocal^\star$. Denote by ${\Ocal^\star}'$ the orientation of $G^\star$ obtained from reversing $e^\star$ in $\Ocal^\star$, then $[(F'^\star)-(F^\star)]+[D_{\Ocal^\star}]=[D_{{\Ocal^\star}'}]$. But it is easy to see ${\Ocal^\star}'$ is the dual orientation of $\Ocal'$, thus $[D_{{\Ocal^\star}'}]$ is the dual element of $[D_{\Ocal'}]$. Summarizing we have the desired identity. \hfill $\Box$\\

The following duality result on Bernardi processes was first observed by Bernardi himself \cite[Section 8.2]{Bernardi_Process}, but for the sake of completeness we include a proof using our language here.

\begin{lemma} \label{Dual_orient}
Let $G$ be a bridgeless planar ribbon graph and let $G^\star$ be a planar dual of $G$. Let $\beta^\star$ be the Bernardi bijection on $G^\star$ corresponding to a dual face $v^\star$. Then for any spanning tree $T$ of $G$ and $e=xy\in T$, $\beta^\star$ orients $e^\star$ the same way as orienting $e$ away from $v$ in $T$ (namely, if $x$ is closer to $v$ than $y$ in $T$, then $e$ is oriented as $\overrightarrow{xy}$ and vice versa).
\end{lemma}

\noindent{\sc Proof:} WLOG $v^\star$ is the unbounded face of $G^\star$. Let $C^\star$ be the fundamental cycle of $e^\star$ with respect to $T^\star$. On the one hand, by the dual of Proposition \ref{Plane_BP_as_OM}, $C^\star$ is oriented counter-clockwise and $e^\star$ will be oriented according to this orientation of $C^\star$ by $\beta^\star$. On the other hand, $C^\star$ bounds a region that does not contain $v^\star$, so orienting $e$ away from $v$ means orienting $e$ to go into the bounded region. Therefore the two methods give the same orientation to $e$ and $e^\star$ in terms of the canonical orientation identification.\hfill $\Box$

\begin{theorem} \label{strong_Bernardi_CD}
Let $G$ be a bridgeless planar ribbon graph and let $G^\star$ be a planar dual of $G$. Let $\beta_F,\beta_{v^\star}$ be Bernardi bijections of $G,G^\star$ correspond to a face $F$ and a dual face $v^\star$, respectively. Then the following diagram commutes:
\[
\begin{CD}
S(G) @>{\beta_{F}-(v)}>> \Pic^{g-1}(G)\\
@VV{\sigma_T}V @VV{\sigma^{g-1}_D=\tau_{G^\star}\circ\sigma_O\circ\tau_G^{-1}}V\\
S(G^\star) @>{\beta_{v^\star}-(F^\star)}>> \Pic^{g^\star-1}(G^\star)\\
\end{CD}
\]
\end{theorem}

\noindent{\sc Proof:} The diagram says that for every spanning tree $T$ of $G$,
$$\tau_{G^\star}^{-1}(\beta_{v^\star}(\sigma_T(T))-(F^\star))=\sigma_O\circ\tau_{G}^{-1}(\beta_{F}(T)-(v)).$$ We claim that the orientation class $\tau_{G}^{-1}(\beta_{F}(T)-(v))$ contains the orientation $\Ocal$, which is obtained by orienting edges in $T^\star$ (identified as edges in $G$) using $\beta_{F}$ and edges in $T$ using $\beta_{v^\star}$ ($\beta_{v^\star}$ orients edges outside $T^\star$ in $G^\star$, which can be identified as edges in $T$). To see this, note first that if $G$ is partially oriented using $\beta_F$, then $\beta_{v^\star}$ orients edges in $T$ away from $v$ by Lemma \ref{Dual_orient}, which gives an extra in-going edge for every vertex of $G$ except $v$, thus $\tau_G([\Ocal])=[\beta_{F}(T)-(v)]$. Dually $\tau_{G^\star}^{-1}(\beta_{v^\star}(\sigma_T(T))-(F))$ contains the orientation $\Ocal'$ of $G^\star$, which is obtained by orienting edges in $T$ using $\beta_{v^\star}$ and edges in $T^\star$ using $\beta_{F}$. Now visibly $[\Ocal']=\sigma_O([\Ocal])$, as $\Ocal,\Ocal'$ are produced from the same procedure up to duality.\hfill $\Box$\\

\noindent{\sc Proof of Theorem \ref{BW_Dual}:} The diagram there factors through the diagram (and its reversal) in Theorem \ref{strong_Bernardi_CD} and the diagram in Proposition \ref{Pic_CD} as follows.\\

\noindent\scaleto{
\minCDarrowwidth20pt
\begin{CD}
\Pic^0(G) \times S(G) @>{{\rm id}\times(\beta_{F}-(v))}>> \Pic^0(G) \times\Pic^{g-1}(G) @>>> \Pic^{g-1}(G) @>{\beta_F^{-1}(\cdot+(v))}>> S(G) \\
@VV{\sigma^0_D \times \sigma_T}V	@VV{\sigma^0_D\times\sigma^{g-1}_D}V	@VV{\sigma^{g-1}_D}V	@VV{\sigma_T}V \\
\Pic^0(G^\star) \times S(G^\star) @>{{\rm id}\times(\beta_{v^\star}-(F^\star))}>> \Pic^0(G^\star) \times\Pic^{g^\star-1}(G^\star)@>>> \Pic^{g^\star-1}(G^\star) @>{\beta_{v^\star}^{-1}(\cdot+(F^\star))}>> S(G^\star)\\
\end{CD}
}{55pt}\\

\hfill $\Box$

\subsection{A Partial Converse for Non-Planar Cases} \label{sec:Non_Planar}

Based on the new results and proofs from the last two sections, one might ask whether Bernardi bijections on non-planar ribbon graphs are also geometric, or at least whether they are cycle orientation maps. One could also ask the opposite question, whether planarity is a necessary condition for a Bernardi bijection to be geometric. In this section, we show that both assertions are incorrect, but the truth is closer to the latter.\\

First we characterize non-planar ribbon graphs by forbidden subdivisions. We say a ribbon graph $\hat{H}$ is a {\em subdivision} of the ribbon graph $H$ if one can obtain $\hat{H}$ from $H$ by inserting degree 2 vertices (equipped with the unique cyclic ordering on two edges) inside the edges of $H$ while keeping the cyclic orderings of edges around the original vertices the same. Also, we say a ribbon graph $H$ is a {\em subgraph} of a ribbon graph $G$ if graph-theoretically $H$ is a subgraph of $G$ and the cyclic ordering of edges around each vertex of $H$ is inherited from the cyclic ordering of $G$. Finally, we say a ribbon graph $G$ contains a ribbon graph $H$ as a subdivision if $\hat{H}$ is a subgraph of $G$ for some subdivision $\hat{H}$ of $H$.

\begin{definition}
The first basic non-planar ribbon graph (BNG I) is a ribbon graph on two vertices $v_1,v_2$ and three edges $e_1,e_2,e_3$ so that the cyclic ordering of the edges around each vertex is $e_1,e_2,e_3$. The second basic non-planar ribbon graph (BNG II) is a ribbon graph on three vertices $u_1,u_2,u_3$, two edges $e'_1,e'_2$ between $u_1,u_2$ and two edges $e'_3,e'_4$ between $u_1,u_3$, where the cyclic ordering of the edges around $u_1$ is $e'_1,e'_3,e'_2,e'_4$.
\end{definition}

\begin{figure}[h!]
\begin{center}
    \includegraphics[scale=0.4]{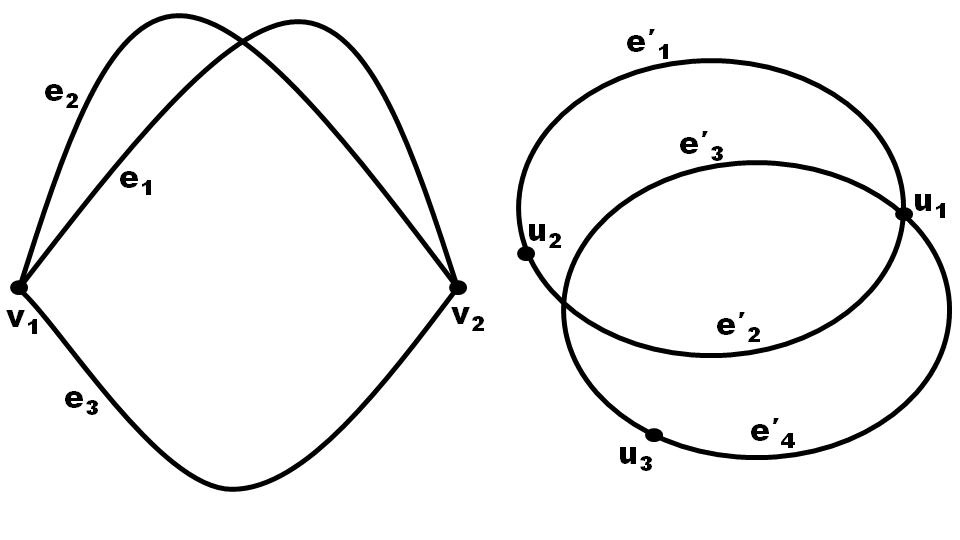}
\end{center}
  \caption{BNG I and BNG II.}
  \label{BNG_diag}
\end{figure}

The next proposition is known within the communities working in structural or topological graph theory. But we could not find an explicit reference in the literature so we include a brief proof here.

\begin{lemma}
Every non-planar ribbon graph contains at least one of the BNGs as a subdivision.
\end{lemma}

\noindent{\sc Proof:} We will demonstrate a way to either find a planar embedding of the ribbon graph $G$, or find a BNG as a subdivision in $G$. First we assume $G$ is 2-connected, and we apply induction via an ear decomposition (a procedure to build any 2-connected graph by starting with a cycle and successively attaching a path (ear) that intersects the current graph exactly at the two distinct endpoints, cf. \cite[Proposition 3.1.3]{Diestel_GT}). The base case is a cycle and is trivial. Suppose a subgraph $G'\subset G$ is embedded in the plane, and let $P=u_0-e_1-\ldots-e_k-u_k$ be an ear to be added. Say in $G'\cup P$ the cyclic ordering of edges around $u_0$ includes $e,e_1,f$ consecutively, then some cycle $C$ containing $e,f$ bounds a face $F$ of $G'$ that $e_1$ is to be embedded in; similarly there is a face $F'$ of $G'$ that $e_k$ is to be embedded in. If $F=F'$ then we can embed $P$ in $F$; otherwise if we let $Q\subset G'$ be a shortest path (possibly trivial) going from $u_k$ to any vertex $v\neq u_0$ on $C$, then $P\cup Q\cup C$ will be a subdivision of BNG I.\\

For the general case, we induct on the number of blocks. Let $v$ be a cut vertex and let $G_1,\ldots,G_k$ be subgraphs corresponding to the components of $T-v$, where $T$ is the block decomposition tree of $G$ (a tree whose vertices are the cut vertices and blocks of $G$, and a cut-vertex $u$ is incident to a block $B$ whenever $u\in B$, cf. \cite[Proposition 3.1.2]{Diestel_GT}). By induction, each $G_i$ can be embedded in the plane, and we may further assume that $v$ is on the boundary of the unbounded face of each embedding if needed. If there exist some interlacing edges $e,e',f,f'$ in the cyclic ordering around $v$ with $e,f\in G_i$ and $e',f'\in G_j$, then by letting $C\subset G_i,C'\subset G_j$ to be cycles containing $e,f$ and $e',f'$, respectively, we have $C\cup C'$ as a subdivision of BNG II in $G$. Otherwise, it can be seen that for all subgraphs $G_i$'s except possibly one (say $G_1$), all edges in $G_i$ incident to $v$ are in some interval $I_i$ of the cyclic ordering around $v$, so we can embed $G_1$ in the plane and then embed other subgraphs one by one according to the cyclic ordering of $I_i$'s. \hfill $\Box$\\

Our partial converse shows that the property ``every Bernardi process is geometric'' is in some sense a property of $G$ as a geometric object, as passing to subdivision does not change the embedding properties of $G$ on surfaces, nor does it change the tropical Jacobian if we keep track of the metric information when we are subdividing edges.\\

We adopt the following conventions: we say an edge $f\neq e_1,e_2$ in a ribbon graph is {\em in between $e_1$ and $e_2$ at $v$} if $v$ is a common endpoint of the three edges, and $f$ goes before $e_2$ in the cyclic ordering of edges around $v$ when listed starting with $e_1$; given a simple path $P$ and vertices $a,b\in V(P)$, we denote by $aPb$ the subpath of $P$ between $a$ and $b$ (inclusive); and given a spanning tree $T$ and a subset of vertices $V'$ such that $T[V']$ is connected, we say a vertex $v$ not in $V'$ is {\em under a vertex $v'\in V'$} if the closest vertex from $V'$ to $v$ in $T$ is $v'$.

\begin{theorem} \label{Non_Planar_BP}
Let $G$ be a non-planar ribbon graph. If $G$ is simple, then there exists some Bernardi bijection on $G$ that is not a cycle orientation map. Otherwise there exists some subdivision $G'$ of $G$ such that some Bernardi bijection on $G'$ is not a cycle orientation map, and $G'$ can be chosen to have at most one more vertex than $G$.
\end{theorem}

\noindent{\sc Proof:} Let $G$ be a non-planar ribbon graph containing a subdivision $P_1\cup P_2\cup P_3$ of BNG I, with $P_1,P_2,P_3$ being internally disjoint paths sharing endpoints $v_1,v_2$; we assume $P_1$ is not an edge in the simple graph case by re-indexing. WLOG we may assume there are no edges between (the last edge of) $P_1$ and (the last edge of) $P_2$ whose endpoints are $v_2$ and some internal vertex of $P_1$ or $P_2$: if there is such an edge $f=v_2t$ with $t\in V(P_i)$, $i=1$ or 2, then we can replace $P_i$ by $v_1P_it\cup\{f\}$, the process will eventually stop because the number of edges between $P_1$ and $P_2$ decreases in every step. Note that in the case of simple graphs, the process will keep at least one internal vertex of $P_1$. Similarly we may assume there are no edges between $v_1$ and $v_2$ that are between $P_1,P_2$ at the two ends, or otherwise we may replace $P_2$ by such edge. By inserting a new vertex on $P_1$ near $v_2$ in the non-simple case if necessary, we may assume $P_1=v_1-\ldots-e_{11}-u-e_{12}-v_2$ is of length at least 2, and there are no edges between $u$ and $v_2$ other than $e_{12}$.\\

Denote by $e_3$ the edge on $P_3$ that is incident to $v_2$, we extend the acyclic subgraph $(P_1-e_{11})\cup P_2\cup (P_3-e_3)$ to a spanning tree $T_1$ of $G$ with the maximum number of vertices under $v_2$ with respect to $V[P_1\cup P_2\cup P_3]$. Note that our assumption means the other endpoint of any non-tree edge incident to $v_2$ in $T_1$ is either from $V[P_1\cup P_2\cup P_3]$ or is a vertex under $v_2$. Set $T_2=T_1-e_{12}+e_{11}$. It is easy to see that the set of non-tree edges incident to $v_2$ in $T_2$ is exactly the set of non-tree edges incident to $v_2$ in $T_1$ plus $e_{12}$, and those common non-tree edges have the same fundamental cycles in $T_1$ and $T_2$. Consider the Bernardi tours of $T_1,T_2\subset G$ starting from $(u,e_{11})$. A routine simulation shows that each non-tree edge that contributes a chip to $v_2$ in the first tour will also contribute a chip to $v_2$ in the second, while $e_{12}$ and $e_3$ will each contribute a chip to $v_2$ in the second tour but not in the first, so $v_2$ received at least two more chips in the second tour than in the first. But in any cycle orientation map, $v_2$ can receive at most one more chip with respect to $T_2$ than $T_1$, so the Bernardi bijection is not a cycle orientation map.\\

The case when $G$ only has subdivisions of BNG II is similar. Let $C_1=u_1-e_1-u-f-\ldots-e_3-u_1$ ($f$ could be equal to $e_3$) and $C_2=u_1-e_2-\ldots-w-e_4-u_1$ be two cycles of $G$ whose union is a subdivision of BNG II, i.e. the two cycles are disjoint except at $u_1$, and the cyclic ordering of edges around $u_1$ includes $e_1,e_2,e_3,e_4$ in order. With a greedy procedure similar to the one in the case of BNG I, we may assume there are no edges between $e_1$ and $e_2$ whose endpoints are $u_1$ and some other vertex of $C_1$ or $C_2$, nor edges between $e_2$ and $e_3$ whose endpoints are $u_1$ and some other vertex of $C_1$.  By inserting a new vertex in $e_1$ near $u_1$ in the non-simple case, we may further assume there are no edges between $u_1$ and $u$ other than $e_1$. Extend $(C_1-f)\cup (C_2-e_4)$ to a spanning tree $T'_1$ of $G$ with the maximum number of vertices under $u_1$ with respect to $V[C_1\cup C_2]$, and set $T'_2=T'_1-e_1+f$. Consider the two Bernardi tours of $T'_1,T'_2$ starting from $(u,f)$. $u_1$ received at least two more chips (from $e_1$ and $e_4$) in the second tour than in the first, a similar reasoning as above shows the Bernardi bijection can not be a cycle orientation map. \hfill $\Box$\\

We end this section by noting that it can be checked that all Bernardi bijections on BNG I are geometric, so working with subdivisions is indeed necessary in Theorem \ref{Non_Planar_BP}.

\section{Conclusion and Open Problems} \label{sec:Problems}

We started with a geometric object rooted in tropical geometry and discovered several classes of bijections between spanning trees and other combinatorial objects which have nice combinatorial and algorithmic properties. The new bijections have surprising connections to other bijections defined in quite different ways. We conclude with a few open questions.\\

To motivate some of the following questions, we list several observations on the polyhedral decomposition without proof.

\begin{proposition} \label{GB_Torsor}
Suppose $g>1$. Let $\beta,\beta '$ be edge ordering maps with data $e_1<e_2<\ldots<e_m$ and $\overline{e}_1<e_2<\ldots<e_h$, respectively. Let $C,C'$ be the cones in the geometric bijection fan corresponding to $\beta,\beta'$. Then $C'$ intersects the antipodal cone of $C$ by some positive dimensional face, which contains the ray spanned by $\pi(e_1)$.
\end{proposition}

\begin{question}
When do two geometric bijections give isomorphic torsors? Is it true that for any full-dimensional cone $C$ of the geometric bijection fan, there exists some cone adjacent to the antipodal cone of $C$ that gives an isomorphic torsor?
\end{question}

\begin{proposition} \label{f_vector}
The number of $i$-dimensional faces in the polyhedral decomposition of $\Pic^g(\Gamma)$ is equal to $\binom{g}{i}$ times the number of spanning trees of $G$.
\end{proposition}

\begin{question}
Can one describe a ``geometric bijection'' between the set of $i$-dimensional faces and the set of $(g-i)$-dimensional faces of the polyhedral decomposition, or at least give a high-level explanation for the equality of cardinalities? (The generic shifting setup does not work in general.)
\end{question}

\begin{question}
Classify all cycle orientation configurations that give bijective cycle orientation maps.
\end{question}

\section{Appendix: A Generic Example of Geometric Cycle Orientation Configuration}

\begin{figure}[h!]
\begin{center}
    \includegraphics[scale=0.4]{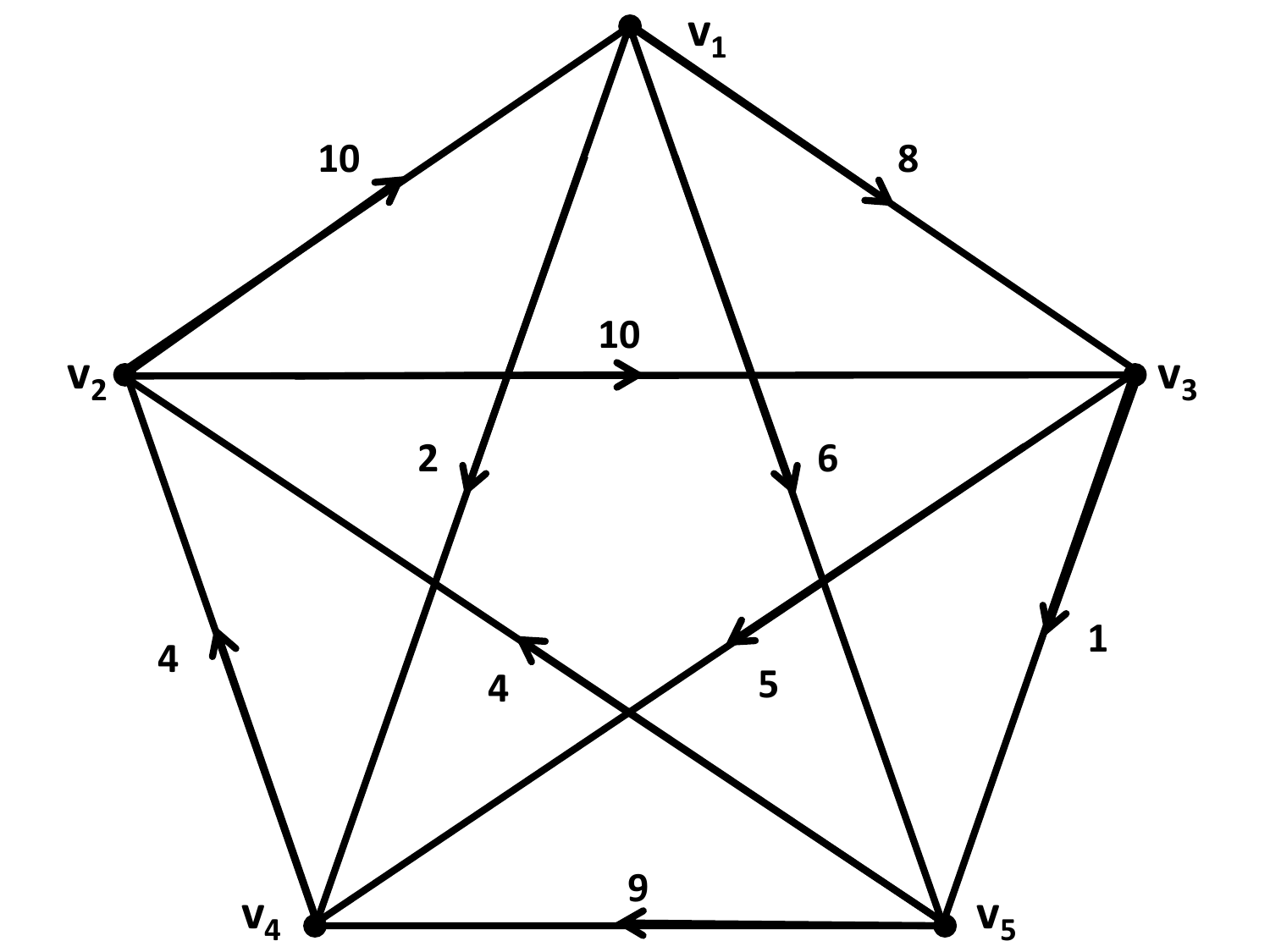}
\end{center}
  \caption{Weights of $E(K_5)$}
  \label{K5_Weight}
\end{figure}

We give an example of cycle orientation configuration that is geometric but not induced by any edge ordering. Consider the weights given to the edges of the complete graph on 5 vertices as in Figure \ref{K5_Weight}. It is routine to check that the weights are generic, i.e. the sum of weights along each cycle is non-zero, hence they induce a geometric cycle orientation configuration. In fact, we can see the weights induce the following directed cycles:\\
$\overrightarrow{v_1v_3v_2}, \overrightarrow{v_1v_4v_2}, \overrightarrow{v_1v_5v_2}, \overrightarrow{v_1v_3v_4}, \overrightarrow{v_1v_3v_5}, \overrightarrow{v_1v_5v_4}, \overrightarrow{v_2v_3v_4}, \overrightarrow{v_2v_5v_4}, \overrightarrow{v_3v_5v_4},\\ \overrightarrow{v_1v_2v_3v_4}, \overrightarrow{v_1v_5v_3v_2},\overrightarrow{v_1v_4v_5v_2},
\overrightarrow{v_1v_4v_2v_3},\overrightarrow{v_1v_5v_2v_4}.$\\

Now notice that each edge appeared in both directions at least once in the above list, so the cycle orientation configuration can not be induced by an ordering of edges, for otherwise the smallest edge will always go the same direction.

School of Mathematics, Georgia Institute of Technology

Atlanta, Georgia 30332-0160, USA

Email address: \url{cyuen7@math.gatech.edu}

\end{document}